\documentclass[12pt,draftcls,onecolumn]{IEEEtran}

\usepackage{amssymb,amsmath}
\usepackage{theorem,bbm}
\usepackage[usenames,dvipsnames]{color}
\usepackage{graphicx,dblfloatfix}
\usepackage{mathrsfs}
\usepackage{subfigure}

\newtheorem{theorem}{Theorem}[section]
\newtheorem{lemma}[theorem]{Lemma}
\newtheorem{remark}[theorem]{Remark}
\newtheorem{corollary}[theorem]{Corollary}
\newtheorem{proposition}[theorem]{Proposition}
\newtheorem{definition}[theorem]{Definition}

\newcommand{\reals}{\mathbb{R}} 
\newcommand{\realsnonnegative}{\mathbb{R}_{\ge 0}}
 
\newcommand{\ball}{\mathbb{B}}

\newcommand{\bd}{\operatorname{bd}}

\newcommand{\diag}{\operatorname{diag}}
\newcommand{\argmin}{\operatorname{argmin}}

\newcommand{\esssup}{\operatorname{esssup}}
\newcommand{\proj}{\operatorname{proj}}
\newcommand{\graph}{\mathcal{G}} 
\newcommand{\primalsol}{\mathcal{X}} 
\newcommand{\pdsol}{\primalsol \times \dualsol} 
\newcommand{\dualsol}{\mathcal{Z}} 
\newcommand{\vertices}{\mathcal{V}}
 
\newcommand{\edges}{\mathcal{E}}

\newcommand{\nonnegint}{\mathbb{Z}_{\ge 0}}
\newcommand{\sharedagent}{\mathbb{S}}
\newcommand{\failededges}{\mathbb{F}}
\newcommand{\sdl}{F^K_{\text{sdl}}}
\newcommand{\sdlw}{f^{w}_{\text{dis}}}

\newcommand{\convexhull}{\operatorname{co}}
\newcommand{\Tdiscon}{T_{\operatorname{disconnected}}^{\max}}
\newcommand{\Tcon}{T_{\operatorname{connected}}^{\min}}

\newcommand{\Vcirc}{V_{\operatorname{euc}}}
\newcommand{\Vcircsep}{V_{\operatorname{euc}}^{\operatorname{rep}}}
\newcommand{\Vclfsep}{V_{\operatorname{CLF}}^{\operatorname{rep}}}
\newcommand{\Vclf}{V_{\operatorname{CLF}}}
\newcommand{\nomflow}{f^{\operatorname{nom}}}
\newcommand{\anomflow}{f^{\operatorname{nom,RCG}}}

\newcommand{\oprocendsymbol}{\hbox{$\bullet$}}
\newcommand{\oprocend}{\relax\ifmmode\else\unskip\hfill\fi\oprocendsymbol}

\newcommand{\until}[1]{\{1,\dots,#1\}}
\newcommand{\MM}{\mathcal{M}}
 
\renewcommand{\epsilon}{\varepsilon} 

\newcommand{\longthmtitle}[1]{\mbox{}\textup{\bf {(#1).}}}
\newcommand{\NormInf}[1]{\|#1\|_{\infty}}
\newcommand{\NormTwo}[1]{\|#1\|}
\newcommand{\setdef}[2]{\{#1 \; |\; #2\}}

\newcommand{\myclearpage}{\clearpage}
\renewcommand{\myclearpage}{}

\begin{document}

\title{Robust distributed linear programming\footnote{Incomplete
    versions of this paper were submitted to the 2013 American Control
    Conference and the 2013 IEEE Conference on Decision and Control.}}





\author{Dean Richert \qquad Jorge Cort\'es \thanks{The authors are
    with the Department of Mechanical and Aerospace Engineering,
    University of California, San Diego, CA 92093, USA, {\tt\small
      \{drichert,cortes\}@ucsd.edu}}}

\maketitle

\begin{abstract}
  This paper presents a robust, distributed algorithm to solve general
  linear programs. The algorithm design builds on the characterization
  of the solutions of the linear program as saddle points of a
  modified Lagrangian function.  We show that the resulting
  continuous-time saddle-point algorithm is provably correct but, in
  general, not distributed because of a global parameter associated
  with the nonsmooth exact penalty function employed to encode the
  inequality constraints of the linear program. This motivates the
  design of a discontinuous saddle-point dynamics that, while enjoying
  the same convergence guarantees, is fully distributed and scalable
  with the dimension of the solution vector.
  We also characterize the robustness against disturbances and link
  failures of the proposed dynamics.  Specifically, we show that it is
  integral-input-to-state stable but not input-to-state stable. The
  latter fact is a consequence of a more general result, that we also
  establish, which states that no algorithmic solution for linear
  programming is input-to-state stable when uncertainty in the problem
  data affects the dynamics as a disturbance. Our results allow us to
  establish the resilience of the proposed distributed dynamics to
  disturbances of finite variation and recurrently disconnected
  communication among the agents.  Simulations in an optimal control
  application illustrate the results.
\end{abstract}

\section{Introduction}\label{sec:introduction}


Linear optimization problems, or simply linear programs, model a broad
array of engineering and economic problems and find numerous
applications in diverse areas such as operations research, network
flow, robust control, microeconomics, and company management. In this
paper, we are interested in both the synthesis of distributed
algorithms that can solve standard form linear programs and the
characterization of their robustness properties. Our interest is
motivated by multi-agent scenarios that give rise to linear programs
with an intrinsic distributed nature.  In such contexts, distributed
approaches have the potential to offer inherent advantages over
centralized solvers. Among these, we highlight the reduction on
communication and computational overhead, the availability of simple
computation tasks that can be performed by inexpensive and
low-performance processors, and the robustness and adaptive behavior
against individual failures.  Here, we consider scenarios where
individual agents interact with their neighbors and are only
responsible for computing their own component of the solution vector
of the linear program. We study the synthesis of provably correct,
distributed algorithms that make the aggregate of the agents' states
converge to a solution of the linear program and are robust to
disturbances and communication link failures.

\emph{Literature review.} Linear programs play an important role in a
wide variety of applications, including perimeter
patrolling~\cite{RA-RC-AC-LS:12}, task
allocation~\cite{DPB:98,MJ-SA-ME:06}, operator
placement~\cite{BWC-AB-FD-BK-KR-GS-RB-MB-FA:12}, process
control~\cite{DRK-JP:64}, routing in communication
networks~\cite{JT-ZH-MJ:07}, and portfolio
optimization~\cite{WFS:67}. This relevance has historically driven the
design of efficient methods to solve linear optimization problems, see
e.g.,~\cite{GBD:63,DB-JNT:97,SB-LV:09}. More recently, the interest on
networked systems and multi-agent coordination has stimulated the
synthesis of distributed strategies to solve linear
programs~\cite{MB-GN-FB-FA:12,GN-FB:11,GY-RS:09} and more general
optimization problems with constraints, see
e.g.,~\cite{AN-AO-PAP:10,MZ-SM:12,JW-NE:11} and references therein.
The aforementioned works build on consensus-based
dynamics~\cite{ROS-JAF-RMM:07,WR-RWB:08,FB-JC-SM:08cor,MM-ME:10}
whereby individual agents agree on the global solution to the
optimization problem. This is a major difference with respect to our
work here, in which each individual agent computes only its own
component of the solution vector by communicating with its neighbors.
This feature makes the messages transmitted over the network
independent of the size of the solution vector, and hence scalable (a
property which would not be shared by a consensus-based distributed
optimization method for the particular class of problems considered
here). Some algorithms that enjoy a similar scalability property exist
in the literature. In particular, the recent work~\cite{RC-GN:13}
introduces a partition-based dual decomposition algorithm for network
optimization.  Other discrete-time algorithms for non-strict convex
problems are proposed in~\cite{DPB-JNT:97,IN-JS:08}, but require at
least one of the exact solutions of a local optimization problem at
each iteration, bounded feasibility sets, or auxiliary variables that
increase the problem dimension. The algorithm in~\cite{NG-JK:98} on
the other hand only achieves convergence to an approximate solution of
the optimization problem. Closer to our approach, although without
equality constraints, the works~\cite{DF-FP:10,KA-LH-HU:58} build on
the saddle-point dynamics of a smooth Lagrangian function to propose
algorithms for linear programming.  The resulting dynamics are
discontinuous because of the projections taken to keep the evolution
within the feasible set. Both works establish convergence in the
primal variables under the assumption that the solution of the linear
program is unique~\cite{KA-LH-HU:58} or that Slater's condition is
satisfied~\cite{DF-FP:10}, but do not characterize the properties of
the final convergence point in the dual variables, which might indeed
not be a solution of the dual problem. We are unaware of works that
explicitly address the problem of studying the robustness of linear
programming algorithms, particularly from a systems and control
perspective. This brings up another point of connection of the present
treatment with the literature, which is the body of work on robustness
of dynamical systems against disturbances. In particular, we explore
the properties of our proposed dynamics with respect to notions such
as robust asymptotic stability~\cite{CC-ART-RG:08}, input-to-state
stability (ISS)~\cite{EDS:89-tac}, and integral input-to-state
stability (iISS)~\cite{DA-EDS-YW:00}.  The term `robust optimization'
often employed in the literature, see e.g.~\cite{DB-DBB-CC:11}, refers
instead to worst-case optimization problems where uncertainty in the
data is explicitly included in the problem formulation. In this
context, `robust' refers to the problem formulation and not to the
actual algorithm employed to solve the optimization.


\emph{Statement of contributions.} We consider standard form linear
programs, which contain both equality and non-negativity constraints
on the decision vector. Our first contribution is an alternative
formulation of the primal-dual solutions of the linear program as
saddle points of a modified Lagrangian function. This function
incorporates an exact nonsmooth penalty function to enforce the
inequality constraints. Our second contribution concerns the design of
a continuous-time dynamics that find the solutions of standard form
linear programs.  Our alternative problem formulation motivates the
study of the saddle-point dynamics (gradient descent in one variable,
gradient ascent in the other) associated with the modified Lagrangian.
It should be noted that, in general, saddle points are only guaranteed
to be stable (and not necessarily asymptotically stable) for the
corresponding saddle-point dynamics. Nevertheless, in our case, we are
able to establish the global asymptotic stability of the (possibly
unbounded) set of primal-dual solutions of the linear program and,
moreover, the pointwise convergence of the trajectories.  Our analysis
relies on the set-valued LaSalle Invariance Principle and, in
particular, a careful use of the properties of weakly and strongly
invariant sets of the saddle-point dynamics.  In general, knowledge of
the global parameter associated with the nonsmooth exact penalty
function employed to encode the inequality constraints is necessary
for the implementation of the saddle-point dynamics. To circumvent
this need, we propose an alternative discontinuous saddle-point
dynamics that does not require such knowledge and is fully distributed
over a multi-agent system in which each individual computes only its
own component of the solution vector. We show that the discontinuous
dynamics share the same convergence properties of the regular
saddle-point dynamics by establishing that, for sufficiently large
values of the global parameter, the trajectories of the former are
also trajectories of the latter. Two key advantages of our methodology
are that it (i) allows us to establish global asymptotic stability of
the discontinuous dynamics without establishing any regularity
conditions on the switching behavior and (ii) sets the stage for the
characterization of novel and relevant algorithm robustness
properties. This latter point bring us to our third contribution,
which pertains the robustness of the discontinuous saddle-point
dynamics against disturbances and link failures. We establish that no
continuous-time algorithm that solves general linear programs can be
input-to-state stable (ISS) when uncertainty in the problem data
affects the dynamics as a disturbance. As our technical approach
shows, this fact is due to the intrinsic properties of the primal-dual
solutions to linear programs.
Nevertheless, when the set of primal-dual solutions is compact, we
show that our discontinuous saddle-point dynamics possesses an
ISS-like property against small constant disturbances and, more
importantly, is integral input-to-state stable (iISS) -- and thus
robust to finite energy disturbances. Our proof method is based on
identifying a suitable iISS Lyapunov function, which we build by
combining the Lyapunov function used in our LaSalle argument and
results from converse Lyapunov theory. We conclude that one cannot
expect better disturbance rejection properties from a linear
programming algorithm than those we establish for our discontinuous
saddle-point dynamics. These results allow us to establish the
robustness of our dynamics against disturbances of finite variation
and communication failures among agents modeled by recurrently
connected graphs. Simulations in an optimal control problem illustrate
the results.

\emph{Organization.} Section~\ref{sec:prelim} introduces basic
preliminaries. Section~\ref{sec:problem_statement} presents the
problem statement.  Section~\ref{sec:DLP} proposes the discontinuous
saddle-point dynamics, establishes its convergence, and discusses its
distributed implementation. Sections~\ref{sec:robustness}
and~\ref{sec:connected} study the algorithm robustness against
disturbances and communication link failures, respectively.
Simulations illustrate our results in
Section~\ref{sec:simulations}. Finally, Section~\ref{sec:conclusions}
summarizes our results and ideas for future work.

\myclearpage

\section{Preliminaries}\label{sec:prelim}

Here, we introduce notation and basic notions on nonsmooth analysis
and dynamical systems. This section may be safely skipped by the
reader who is familiar with the notions reviewed here.

\subsection{Notation and basic notions}

The set of real numbers is $\reals$.  For $x \in \reals^n$,
$x \ge 0$ (resp. $x > 0$) means that all components of $x$ are
nonnegative (resp. positive). For $x \in \reals^n$, we define $\max \{
0, x \} = (\max\{0,x_1\},\dots,\max\{0,x_n\}) \in \reals_{\ge 0}^n$.
We let $\mathbbm{1}_n \in \reals^n$ denote the vector of ones.  We use
$\NormTwo{\cdot}$ and $\NormInf{\cdot}$ to denote the $2$- and
$\infty$-norms in $\reals^n$. The Euclidean distance from a point $x
\in \reals^n$ to a set $A \subset\reals^n$ is denoted by
$\NormTwo{\cdot}_{A}$. The set $\ball(x,\delta) \subset \reals^n$ is
the open ball centered at $x \in \reals^n$ with radius $\delta > 0$.
The set $A \subset \reals^n$ is convex if it fully contains the
segment connecting any two points in $A$.

A function $V: \reals^n \rightarrow \reals$ is positive definite with
respect to $A \subset \reals^n$ if (i) $V(x) = 0$ for all $x \in A$
and $V(x) > 0$ for all $x \notin A$. If $A = \{0\}$, we refer to $V$
as positive definite. $V:\reals^n \rightarrow \reals$ is radially
unbounded with respect to $A$ if $V(x) \rightarrow \infty$ when
$\NormTwo{x}_A \rightarrow \infty$. If $A = \{0\}$, we refer to $V$ as
radially unbounded.  A function $V$ is proper with respect to $A$ if
it is both positive definite and radially unbounded with respect to
$A$. A set-valued map $F:\reals^n \rightrightarrows \reals^n$ maps
elements in $\reals^n$ to subsets of $\reals^n$.  A function $V:X
\rightarrow \reals$ defined on the convex set $X \subset \reals^n$ is
convex if $V(kx + (1-k)y) \le kV(x) + (1-k)V(y)$ for all $x,y \in X$
and $k \in [0,1]$. $V$ is concave iff $-V$ is convex. Given $\rho \in
\reals$, we define $V^{-1}(\le \rho) = \setdef{x \in X}{V(x)\le
  \rho}$. The function $L:X \times Y \rightarrow \reals$ defined on
the convex set $X \times Y \subset \reals^n \times \reals^m$ is
convex-concave if it is convex on its first argument and concave on
its second. A point $(\bar{x},\bar{y}) \in X \times Y$ is a saddle
point of $L$ if $L(x,\bar{y}) \ge L(\bar{x},\bar{y}) \ge L(\bar{x},y)$
for all $(x,y) \in X \times Y$.

The notion of comparison function is useful to formalize stability
properties.  The class of $\mathcal{K}$ functions is composed by
functions of the form $[0,\infty) \rightarrow [0,\infty)$ that are
continuous, zero at zero, and strictly increasing. The subset of class
$\mathcal{K}$ functions that are unbounded are called
class~$\mathcal{K}_{\infty}$. A class $\mathcal{KL}$ function
$[0,\infty) \times [0,\infty) \rightarrow [0,\infty)$ is class
$\mathcal{K}$ in its first argument and continuous, decreasing, and      
converging to zero in its second argument.

An undirected graph is a pair $\graph = (\vertices,\edges)$, where
$\vertices = \{1, \dots, n\}$ is a set of vertices and $\edges
\subseteq \vertices \times \vertices$ is a set of edges. Given a
matrix $A \in \reals^{m \times n}$, we call a graph \emph{connected
  with respect to $A$} if for each $\ell \in \until{m}$ such that
$a_{\ell,i} \not= 0 \neq a_{\ell,j}$, it holds that $(i,j ) \in
\edges$.

\subsection{Nonsmooth analysis}

Here we review some basic notions from nonsmooth analysis
following~\cite{FHC:83}. A function $V:\reals^n \rightarrow \reals$ is
locally Lipschitz at $x \in \reals^n$ if there exist $\delta_x > 0$
and $L_x \ge 0$ such that $|V(y_1) - V(y_2)| \le L_x
\NormTwo{y_1-y_2}$ for $y_1,y_2 \in \ball(x,\delta_x)$. If $V$ is
locally Lipschitz at all $x \in \reals^n$, we refer to $V$ as locally
Lipschitz. If $V$ is convex, then it is locally Lipschitz. A locally
Lipschitz function is differentiable almost everywhere. Let
$\Omega_V\subset \reals^n$ be then the set of points where $V$ is not
differentiable. The generalized gradient of a locally Lipschitz
function $V$ at $x \in \reals^n$ is
\begin{align*}
  \partial V(x) = \convexhull \Big\{ \lim_{i \rightarrow \infty}
  \nabla V(x_i) : x_i \rightarrow x, x_i \notin S \cup \Omega_V
  \Big\},
\end{align*}
where $\convexhull \{ \cdot \}$ denotes the convex hull and $S \subset
\reals^n$ is any set with zero Lebesgue measure.  A critical point $x
\in \reals^n$ of $V$ satisfies $0 \in \partial V(x)$.  For a convex
function $V$, the first-order condition of convexity states that $V(y)
\ge V(x) + (y-x)^Tg$ for all $g \in \partial V(x)$ and $x,y \in
\reals^n$. For $V: \reals^n \times \reals^n \rightarrow \reals$ and
$(x,y) \in \reals^n \times \reals^n$, we use $\partial_xV(x,y)$ and
$\partial_yV(x,y)$ to denote the generalized gradients of the maps $x'
\mapsto V(x',y)$ at $x$ and $y' \mapsto V(x,y')$ at $y$, respectively.

A set-valued map $F:X \subset \reals^n
\rightrightarrows \reals^n$ is upper semi-continuous if for all $x \in
X$ and $\epsilon \in (0,\infty)$ there exists $\delta_x \in
(0,\infty)$ such that $F(y) \subseteq F(x) + \ball(0,\epsilon)$ for
all $y \in \ball(x,\delta_x)$. Conversely, $F$ is lower
semi-continuous if for all $x \in X$, $\epsilon \in (0,\infty)$, and
any open set $A$ intersecting $F(x)$ there exists a $\delta \in
(0,\infty)$ such that $F(y)$ intersects $A$ for all $y \in
\ball(x,\delta)$. If $F$ is both upper and lower semi-continuous then
it is continuous. Also, $F$ is locally bounded if for every $x \in X$
there exist $\epsilon \in (0,\infty)$ and $M > 0$ such that
$\NormTwo{z} \le M$ for all $z \in F(y)$ and all $y \in
\ball(x,\epsilon)$.

\begin{lemma}\longthmtitle{Properties of the generalized
    gradient}\label{lem:gen}
  If~$V:\reals^n \rightarrow \reals$ is locally Lipschitz at $x \in
  \reals^n$, then $\partial V(x)$ is nonempty, convex, and
  compact. Moreover, $x \mapsto \partial V(x)$ is locally bounded and
  upper semi-continuous.
\end{lemma}

\subsection{Set-valued dynamical systems}

Our exposition on basic concepts for set-valued dynamical systems
follows~\cite{JC:08-csm-yo}. A time-invariant set-valued dynamical
system is represented by the differential inclusion
\begin{align}\label{eq:set_val_dyn}
  \dot{x} \in F(x), 
\end{align}
where $t \in \reals_{\ge 0}$ and $F : \reals^n \rightrightarrows
\reals^n$ is a set valued map.  If $F$ is locally bounded, upper
semi-continuous and takes nonempty, convex, and compact values, then
from any initial condition in~$\reals^n$, there exists an absolutely
continuous curve $x:\reals_{\ge 0} \rightarrow \reals^n$, called
solution, satisfying~\eqref{eq:set_val_dyn} almost everywhere.  The
solution is maximal if it cannot be extended forward in time. The set
of equilibria of $F$ is defined as $\setdef{ x \in \reals^n}{0 \in
  F(x) }$. A set $\mathcal{M}$ is strongly (resp. weakly) invariant
with respect to~\eqref{eq:set_val_dyn} if, for each $x_0 \in
\mathcal{M}$, $\mathcal{M}$ contains all (resp. at least one) maximal
solution(s) of~\eqref{eq:set_val_dyn} with initial condition $x_0$.
The set-valued Lie derivative of a differentiable function $V:\reals^n
\rightarrow \reals$ along the trajectories of~\eqref{eq:set_val_dyn}
is defined as
\begin{align*}
  \mathcal{L}_F V(x) = \{ \nabla V(x)^T v : v \in F(x) \}.
\end{align*}
The following result helps establish the asymptotic convergence
properties of~\eqref{eq:set_val_dyn}.

\begin{theorem}\longthmtitle{Set-valued LaSalle Invariance
    Principle}\label{th:lasalle}
  Let $X \subset \reals^n$ be compact and strongly invariant with
  respect to~\eqref{eq:set_val_dyn}.  Assume $V:\reals^n \rightarrow
  \reals$ is differentiable
  and $F$ is locally bounded, upper semi-continuous and takes
  nonempty, convex, and compact values. If $\mathcal{L}_FV(x) \subset
  (-\infty,0]$ for all $x \in X$, then any solution
  of~\eqref{eq:set_val_dyn} starting in $X$ converges to the largest
  weakly invariant set $\mathcal{M}$ contained in $\overline{\{ x \in
    X : 0 \in \mathcal{L}_FV(x)\}}$.
\end{theorem}

Differential inclusions are specially useful to handle differential
equations with discontinuities. Specifically, let $f:X \subset
\reals^n \rightarrow \reals^n$ be a piecewise continuous vector field
and consider
\begin{align}\label{eq:ode-discontinuous}
  \dot x = f(x) .
\end{align}
The classical notion of solution is not applicable
to~\eqref{eq:ode-discontinuous} because of the
discontinuities. Instead, consider the Filippov set-valued map
associated to $f$, defined by
$ \mathcal{F}[f](x) := \overline{\convexhull} \Big \{ \lim_{i
  \rightarrow \infty} f(x_i) : x_i \rightarrow x, x_i \notin \Omega_f
\Big \}$,
where $\overline{\convexhull} \{ \cdot \}$ denotes the closed convex
hull and $\Omega_f$ are the points where $f$ is discontinuous.  One can
show that the set-valued map $\mathcal{F}[f]$ is locally bounded,
upper semi-continuous and takes nonempty, convex, and compact values,
and hence solutions exist to
\begin{align}\label{eq:odi-Filippov}
  \dot{x} \in \mathcal{F}[f](x)
\end{align} 
starting from any initial condition.  The solutions
of~\eqref{eq:ode-discontinuous} in the sense of Filippov are, by
definition, the solutions of the differential
inclusion~\eqref{eq:odi-Filippov}.



\myclearpage
\section{Problem statement and equivalent
  formulation} \label{sec:problem_statement} 

This section introduces standard form linear programs and describes an
alternative formulation that is useful later in fulfilling our main
objective, which is the design of robust, distributed algorithms to
solve them.  Consider the following standard form linear program,
\begin{subequations}\label{eq:standard_form}
  \begin{alignat}{2}
    &\min && \quad c^Tx 
    \\
    &\hspace{1.5mm}\text{s.t.} && \quad Ax = b,
    \quad x \ge 0, \label{eq:con_2}
  \end{alignat}
\end{subequations}
where $x, c \in \reals^n$, $A \in \reals^{m \times n}$, and $b \in
\reals^m$.  We only consider feasible linear programs with finite
optimal value. The set of solutions to~\eqref{eq:standard_form} is
denoted by $\primalsol \subset \reals^n$.  The dual formulation is
\begin{subequations}\label{eq:dual}
  \begin{alignat}{2}
    &\max && \quad -b^Tz \\ &\hspace{1.5mm}\text{s.t.} && \quad A^Tz +
    c \ge 0. \label{eq:dual_con}
  \end{alignat}
\end{subequations}
The set of solutions to~\eqref{eq:dual} is denoted by $\dualsol
\subset \reals^m$.  We use $x_*$ and $z_*$ to denote a solution
of~\eqref{eq:standard_form} and~\eqref{eq:dual}, respectively. The
following result is a fundamental relationship between primal and dual
solutions of linear programs and can be found in many optimization
references, see e.g.,~\cite{SB-LV:09}.

\begin{theorem}\longthmtitle{Complementary slackness and strong
    duality}\label{th:comp}
  Suppose that $x \in \reals^n$ is feasible
  for~\eqref{eq:standard_form} and $z \in \reals^m$ is feasible
  for~\eqref{eq:dual}. Then $x$ is a solution
  to~\eqref{eq:standard_form} and $z$ is a solution to~\eqref{eq:dual}
  if and only if $(A^Tz+c)^Tx = 0$. In compact form,
  \begin{align}\label{eq:comp}
    \pdsol = \setdef{(x,z) \in \reals^n \times \reals^m}{Ax = b, \; x
      \ge 0, \; A^Tz + c \ge 0, \; (A^Tz+c)^Tx = 0}.
  \end{align}
  Moreover, for any $(x_*,z_*) \in \pdsol$, it holds that $c^Tx_* =
  -b^Tz_*$.
\end{theorem}

The equality $(A^Tz+c)^Tx = 0$ is called the \emph{complementary
  slackness} condition whereas the property that $c^Tx_* = -b^Tz_*$ is
called \emph{strong duality}. One remarkable consequence of
Theorem~\ref{th:comp} is that the set on the right-hand side
of~\eqref{eq:comp} is convex (because $\primalsol \times \dualsol$ is
convex). This fact is not obvious since the complementary slackness
condition is not affine in the variables $x$ and $z$. This observation
will allow us to use a simplified version of Danskin's Theorem (see
Lemma~\ref{th:danskin}) in the proof of a key result of
Section~\ref{sec:robustness}.  The next result establishes the
connection between the solutions of~\eqref{eq:standard_form}
and~\eqref{eq:dual} and the saddle points of a modified Lagrangian
function. Its proof can be deduced from results on penalty functions
that appear in optimization, see e.g.~\cite{DPB-AN-AEO:03}, but we
include it here for completeness and consistency of the presentation.

\begin{proposition}\longthmtitle{Solutions of linear program as saddle
    points}\label{prop:saddle}
  For $K \ge 0$, let $L^K : \reals^{n} \times \reals^m \rightarrow
  \reals$ be defined by
  \begin{align}\label{eq:Lagrangian}
    L^K(x,z) &= c^Tx + \frac{1}{2}(Ax-b)^T(Ax-b) + z^T(Ax-b) + K
    \mathbbm{1}^T_n \max \{0,-x \}.
  \end{align}
  Then, $L^K$ is convex in $x$ and concave (in fact, linear)
  in~$z$. Moreover,
  \begin{enumerate}
  \item if $x_* \in \reals^n$ is a solution
    of~\eqref{eq:standard_form} and $z_*\in \reals^m$ is a solution
    of~\eqref{eq:dual}, then the point $(x_*,z_*)$ is a saddle point
    of $L^K$ for any $K \ge \NormInf{A^Tz_* + c}$,
  \item if $(\bar{x},\bar{z}) \in \reals^n \times \reals^m$ is a saddle
    point of $L^K$ with $K > \NormInf{A^Tz_* + c}$ for some $z_*
    \in \reals^m$ solution of~\eqref{eq:dual}, then $\bar{x} \in
    \reals^n$ is a solution of~\eqref{eq:standard_form}.
  \end{enumerate}
\end{proposition}
\begin{IEEEproof}
  One can readily see from~\eqref{eq:Lagrangian} that $L^K$ is a
  convex-concave function. Let $x_*$ be a solution
  of~\eqref{eq:standard_form} and let $z_*$ be a solution
  of~\eqref{eq:dual}. To show (i), using the characterization of
  $\pdsol$ described in Theorem~\ref{th:comp} and the fact that $K \ge
  \NormInf{A^Tz_* + c}$, we can write for any $x \in \reals^n$,
  \begin{align*}
    L^K(x,z_*) &= c^Tx + (Ax-b)^T(Ax-b) + z_*^T(Ax-b) + K
    \mathbbm{1}^T_n \max \{0,-x \},
    \\
    &\ge c^Tx + z_*^T(Ax-b) +
    (A^Tz_* + c)^T \max \{0,-x \},
    \\
    &\ge c^Tx + z_*^T(Ax-b) - (A^Tz_* + c)^T x,
    \\
    &= c^Tx + z_*^TA(x-x_*) - (A^Tz_* + c)^T (x-x_*),
    \\
    &= c^Tx - c^T(x-x_*) = c^T x_* = L^K(x_*,z_*).
  \end{align*}
  The fact that $L^K(x_*,z) = L^K(x_*,z_*)$ for any $z \in \reals^m$ is
  immediate. These two facts together imply that $(x_*,z_*)$ is a
  saddle point of $L^K$.

  We prove (ii) by contradiction. Let $(\bar{x},\bar{z})$ be a saddle
  point of $L^K$ with $K > \NormInf{A^Tz_* + c}$ for some $z_* \in
  \dualsol$, but suppose $\bar{x}$ is not a solution
  of~\eqref{eq:standard_form}.  Let $x_* \in \primalsol$.  Since for
  fixed $x$, $z \mapsto L^K(x,z)$ is concave and differentiable, a
  necessary condition for $(\bar{x},\bar{z})$ to be a saddle point of
  $L^K$ is that $A\bar{x} - b = 0$. Using this fact, $L^K(x_*,\bar{z}) \ge
  L^K(\bar{x},\bar{z})$ can be expressed as
  \begin{align}\label{eq:contra_1}
    c^Tx_* &\ge c^T\bar{x} + K \mathbbm{1}^T_n
    \max\{0,-\bar{x}\}.
  \end{align}
  Now, if $\bar{x} \ge 0$, then $ c^Tx_* \ge c^T\bar{x}$, and thus
  $\bar{x}$ would be a solution of~\eqref{eq:standard_form}.  If,
  instead, $\bar{x} \not \ge 0$,
  \begin{align*}
    c^T\bar{x} &= c^Tx_* + c^T(\bar{x}-x_*),
    \\
    &= c^Tx_* - z_*^TA(\bar{x}-x_*) + (A^Tz_*+c)^T(\bar{x}-x_*),
    \\
    & = c^Tx_* - z_*^T(A\bar{x}-b) + (A^Tz_*+c)^T\bar{x},
    \\
    & > c^Tx_* - K \mathbbm{1}^T_n \max\{0,-\bar{x}\},
  \end{align*}
  which contradicts~\eqref{eq:contra_1}, concluding the proof.
\end{IEEEproof}

The relevance of Proposition~\ref{prop:saddle} is two-fold. On the one
hand, it justifies searching for the saddle points of~$L^K$ instead of
directly solving the constrained optimization
problem~\eqref{eq:standard_form}.  On the other hand, given that $L^K$
is convex-concave, a natural approach to find the saddle points is via
the associated saddle-point dynamics. However, for an arbitrary
function, such dynamics is known to render saddle points only stable,
not asymptotically stable (in fact, the saddle-point dynamics derived
using the standard Lagrangian for a linear program does not converge
to a solution of the linear program, see
e.g.,~\cite{RD-PAS-RS:58,KA-LH-HU:58}). Interestingly~\cite{KA-LH-HU:58}, the
convergence properties of saddle-point dynamics can be improved using
penalty functions associated with the constraints to augment the cost
function.  In our case, we augment the linear cost function $c^Tx$
with a quadratic penalty for the equality constraints and a nonsmooth
penalty function for the inequality constraints. This results in the
nonlinear optimization problem,
\begin{align*}
  \min_{Ax=b}c^Tx + \|Ax-b\|^2 + K \mathbbm{1}^T_n\max\{0,-x\},
\end{align*}
whose standard Lagrangian is equivalent to $L^K$. We use the nonsmooth
penalty function to ensure that there is an \emph{exact} equivalence
between saddle points of $L^K$ and the solutions
of~\eqref{eq:standard_form}. Instead, the use of smooth penalty
functions such as the logarithmic barrier function used
in~\cite{JW-NE:11}, results only in approximate solutions. In the next
section, we show that indeed the saddle-point dynamics of $L^K$
asymptotically converges to saddle points.

\begin{remark}\longthmtitle{Bounds on the parameter $K$}
  {\rm It is worth noticing that the lower bounds on $K$ in
    Proposition~\ref{prop:saddle} are characterized by certain dual
    solutions, which are unknown a priori.  Nevertheless, our
    discussion later shows that this problem can be circumvented and
    that knowledge of such bounds is not necessary for the design of
    robust, distributed algorithms that solve linear programs.  }
  \oprocend
\end{remark}

\myclearpage
\section{Saddle-point dynamics for distributed linear
  programming}\label{sec:DLP}

In this section, we design a continuous-time algorithm to find the
solutions of~\eqref{eq:standard_form} and discuss its distributed
implementation in a multi-agent system.  We further build on the
elements of analysis introduced here to characterize the robustness
properties of linear programming dynamics in the forthcoming sections.
Building on the result in Proposition~\ref{prop:saddle}, we consider
the saddle-point dynamics (gradient descent in one argument, gradient
ascent in the other) of the modified Lagrangian $L^K$. Our
presentation proceeds by characterizing the properties of this
dynamics and observing its limitations, leading up to the main
contribution, which is the introduction of a discontinuous
saddle-point dynamics amenable to distributed implementation.

The nonsmooth character of $L^K$ means that its saddle-point dynamics
takes the form of the following differential inclusion,
\begin{subequations}\label{eq:flow}
  \begin{align}
    \dot{x} + c + A^T(z + Ax - b) & \in - K \partial \max \{ 0, -x
    \}, \label{eq:flow_x}
    \\
    \dot{z} &= Ax - b. \label{eq:flow_z}
  \end{align}
\end{subequations}
For notational convenience, we use $\sdl: \reals^n \times \reals^m
\rightrightarrows \reals^n \times \reals^m$ to denote the set-valued
vector field which defines the differential inclusion~\eqref{eq:flow}.
The following result characterizes the asymptotic convergence
of~\eqref{eq:flow} to the set of solutions
to~\eqref{eq:standard_form}-\eqref{eq:dual}.

\begin{theorem}\longthmtitle{Asymptotic convergence to the primal-dual
    solution set}\label{th:grad_flow}
  Let $(x_*,z_*) \in \primalsol \times \dualsol$ and define $V :
  \reals^n \times \reals^m \rightarrow \reals_{\ge 0}$ as
  \begin{align*}
    V(x,z) = \frac{1}{2}(x-x_*)^T(x-x_*) + \frac{1}{2}(z-z_*)^T(z-z_*).
  \end{align*}
  For $\infty > K \ge \NormInf{A^Tz_* + c}$, it holds that
  $\mathcal{L}_{\sdl}V(x,z) \subset (-\infty,0]$ for all $(x,z) \in
  \reals^n \times \reals^m$ and any trajectory $t \mapsto (x(t),z(t))$
  of~\eqref{eq:flow} converges asymptotically to the set $\pdsol$.
\end{theorem}
\begin{IEEEproof}
  Our proof strategy is based on verifying the hypotheses of the
  LaSalle Invariance Principle, cf. Theorem~\ref{th:lasalle}, and
  identifying the set of primal-dual solutions as the corresponding
  largest weakly invariant set.  First, note that Lemma~\ref{lem:gen}
  implies that $\sdl$ is locally bounded, upper semi-continuous and
  takes nonempty, convex, and compact values. By
  Proposition~\ref{prop:saddle}(i), $(x_*,z_*)$ is a saddle point
  of~$L^K$ when $K \ge \NormInf{A^Tz_* + c}$. Consider the quadratic
  function $V$ defined in the theorem statement, which is continuously
  differentiable and radially unbounded. Let $a \in \mathcal{L}_{\sdl}
  V(x,z)$. By definition, there exists $v = \!  (- c - A^T(z + Ax-b) -
  g_x,Ax-b) \in \sdl(x,z)$, with $g_x \in K
  \partial \max\{0,-x\}$, such that
  \begin{align}\label{eq:a}
    a = v^T \nabla V(x,z) = (x \! - \! x_*)^T(- \!  c \!  - \!  A^T(z
    \! + \! Ax \! - \!  b) \! - \! g_x) \! + \!  (z-z_*)(Ax-b).
  \end{align}
  Since $L^K$ is convex in its first argument, and $c + A^T(z + Ax-b) +
  g_x \in \partial_x L^K(x,z)$, using the first-order condition of
  convexity, we have
  \begin{align*}
    L^K(x,z) \le L^K(x_*,z) \! + \! (x \! - \! x_*)^T \! \big( c \!  + \!
    A^T(z \! + \! Ax \! - \! b) \! + \! g_x \big).
  \end{align*}
  Since $L^K$ is linear in $z$, we have $L^K(x,z) = L^K(x,z_*) +
  (z-z_*)^T(Ax - b)$. Using these facts in~\eqref{eq:a}, we get
  \begin{align*}
    a \le L^K(x_*,z) - L^K(x,z_*) = L^K(x_*,z) - L^K(x_*,z_*) +
    L^K(x_*,z_*) -L^K(x,z_*) \le 0 ,
  \end{align*}
  since $(x_*,z_*)$ is a saddle point of $L^K$. Since $a$ is
  arbitrary, we deduce that $\mathcal{L}_{\sdl} V(x,z) \subset
  (-\infty,0]$. For any given $\rho \ge 0$, this implies that the
  sublevel set $V^{-1}(\le \rho)$ is strongly invariant with respect
  to~\eqref{eq:flow}. Since $V$ is radially unbounded, $V^{-1}(\le
  \rho)$ is also compact.  The conditions of Theorem~\ref{th:lasalle}
  are then satisfied with $X=V^{-1}(\le \rho)$, and therefore any
  trajectory of~\eqref{eq:flow} starting in $V^{-1}(\le \rho)$
  converges to the largest weakly invariant set $\mathcal{M}$ in $\{
  (x,z) \in V^{-1}(\le \rho) : 0 \in \mathcal{L}_{\sdl} V(x,z) \}$
  (note that for any initial condition $(x_0,z_0)$ one can choose a
  $\rho$ such that $(x_0,z_0) \in V^{-1}(\le \rho)$). This set is
  closed, which can be justified as follows.  Since $\sdl$ is upper
  semi-continuous and $V$ is continuously differentiable, the map
  $(x,z) \mapsto \mathcal{L}_{\sdl} V(x,z)$ is also upper
  semi-continuous. Closedness then follows from~\cite[Convergence
  Theorem]{JPA-AC:84}. We now show that $\mathcal{M} \subseteq
  \pdsol$. To start, take $(x',z') \in \mathcal{M}$. Then $
  L^K(x_*,z_*) - L^K(x',z_*) = 0$, which implies
  \begin{align}
    \tilde{L}^K(x',z_*) - (Ax'-b)^T(Ax'-b) = 0, \label{eq:lasalle_arg}
  \end{align}
  where $ \tilde{L}^K(x',z_*) = c^Tx_* - c^Tx' - z_*^T(Ax'-b) - K
  \mathbbm{1}^T_n \max\{0,-x' \}$. Using strong duality, the
  expression of~$\tilde{L}^K$ can be simplified
  to~$\tilde{L}^K(x',z_*) = -(A^Tz_* + c)^T x' - K \mathbbm{1}^T_n
  \max\{0,-x' \}$.  In addition, $A^Tz_* + c \ge 0$ by dual
  feasibility. Thus, when $K \ge \NormInf{A^Tz_* + c}$, we have
  $\tilde{L}^K(x,z_*) \le 0$ for all $(x,z) \in V^{-1}(\le
  \rho)$. This implies that $(Ax'-b)^T(Ax'-b) = 0$
  for~\eqref{eq:lasalle_arg} to be true, which further implies that
  $Ax' - b = 0$. Moreover, from the definition of $\tilde{L}^K$ and
  the bound on $K$, one can see that if $x' \not \ge 0$, then
  $\tilde{L}^K(x',z_*) < 0$.  Therefore, for~\eqref{eq:lasalle_arg} to
  be true, it must be that $x' \ge 0$.  Finally,
  from~\eqref{eq:lasalle_arg}, we get that $\tilde{L}^K(x',z_*) =
  c^Tx_*-c^Tx' = 0$.  In summary, if $(x',z') \in \mathcal{M}$ then
  $c^Tx_* = c^Tx'$, $Ax' - b = 0$, and $x' \ge 0$. Therefore, $x'$ is
  a solution of~\eqref{eq:standard_form}.
  Now, we show that $z'$ is a solution of~\eqref{eq:dual}. Because
  $\MM$ is weakly invariant, there exists a trajectory starting from
  $(x',z')$ that remains in $\MM$. The fact that $Ax' = b$ implies
  that $\dot{z} = 0$, and hence $z(t) = z'$ is constant.  For any
  given $i \in \{1,\dots,n\}$, we consider the cases (i) $x'_{i} > 0$
  and (ii) $x'_{i} = 0$. In case~(i), the dynamics of the $i$th
  component of $x$ is $\dot{x}_{i} = -(c + A^Tz')_i$ where $(c +
  A^Tz')_i$ is constant. It cannot be that $ -(c + A^Tz')_i > 0$
  because this would contradict the fact that $t \mapsto x_{i}(t)$ is
  bounded. Therefore, $(c + A^Tz')_i \ge 0$. If $\dot{x}_i = -(c +
  A^Tz')_i < 0$, then $x_{i}(t)$ will eventually become zero, which we
  consider in case (ii). In fact, since the solution remains in $\MM$,
  without loss of generality, we can assume that $(x',z')$ is such
  that either $x'_{i} > 0$ and $(c + A^Tz')_i = 0$ or $x'_{i} = 0$ for
  each $i \in \{1,\dots,n\}$. Consider now case (ii). Since $x_{i}(t)$
  must remain non-negative in $\mathcal{M}$, it must be that
  $\dot{x}_{i}(t) \ge 0$ when $x_{i}(t) = 0$. That is, in
  $\mathcal{M}$, we have $\dot{x}_{i}(t) \ge 0$ when $x_{i}(t) = 0$
  and $\dot{x}_{i}(t) \le 0$ when $x_{i}(t) > 0$. Therefore, for any
  trajectory $t \mapsto x_i(t)$ in $\mathcal{M}$ starting at $x'_{i} =
  0$, the unique Filippov solution is that $x_{i}(t) = 0$ for all $t
  \ge 0$. As a consequence, $(c+A^Tz')_i \in [0,K]$ if $x'_{i} =
  0$. To summarize cases (i) and (ii), we have
  \begin{itemize}
  \item $Ax' = b$ and $x' \ge 0$ (primal feasibility),
  \item $A^Tz' + c \ge 0$ (dual feasibility),
  \item $(A^Tz' + c)_i = 0$ if $x'_{i} > 0$ and $x'_{i} = 0$ if $(A^Tz'
    + c)_i > 0$ (complementary slackness),
  \end{itemize}
  which is sufficient to show that $z \in \dualsol$
  (cf. Theorem~\ref{th:comp}). Hence $\mathcal{M} \subseteq
  \pdsol$. Since the trajectories of~\eqref{eq:flow} converge to
  $\mathcal{M}$, this completes the proof.
\end{IEEEproof}

Using a slightly more complicated lower bound on the parameter $K$, we
are able to show point-wise convergence of the saddle-point
dynamics. We state this result next.

\begin{corollary}\longthmtitle{Point-wise convergence of saddle-point
    dynamics}\label{cor:point}
  Let $\rho > 0$. Then, with the notation of
  Theorem~\ref{th:grad_flow}, for
  \begin{align}\label{eq:gamma_bound}
    \infty > K \ge \max_{(x,z) \in (\primalsol \times
      \dualsol) \cap V^{-1}(\le \rho)}\NormInf{A^Tz + c},
  \end{align}
  it holds that any trajectory $t \mapsto (x(t),z(t))$
  of~\eqref{eq:flow} starting in $V^{-1}(\le \rho)$ converges
  asymptotically to a point in $\pdsol$.
\end{corollary}
\begin{IEEEproof}
  If $K$ satisfies~\eqref{eq:gamma_bound}, then in particular $K \ge
  \NormInf{A^Tz_*+c}$. Thus, $V^{-1}(\le \rho)$ is strongly invariant
  under~\eqref{eq:flow} since $\mathcal{L}_{\sdl}V(x,z) \subset
  (-\infty,0]$ for all $(x,z) \in V^{-1}(\le \rho)$ (cf.
  Theorem~\ref{th:grad_flow}). Also, $V^{-1}(\le \rho)$ is bounded
  because $V$ is quadratic. Therefore, by the Bolzano-Weierstrass
  theorem~\cite[Theorem 3.6]{WR:53}, there exists a subsequence
  $(x(t_k),z(t_k)) \in V^{-1}(\le \hspace*{-2pt} \rho)$ that converges
  to a point $(\tilde{x},\tilde{z}) \in (\pdsol) \cap V^{-1}(\le
  \hspace*{-2pt} \rho)$.  Given $\epsilon>0$, let $k^*$ be such that
  $\NormTwo{(x(t_{k^*}),z(t_{k^*})) - (\tilde{x},\tilde{z})} \le
  \epsilon$. Consider the function $\tilde{V}(x,z) =
  \frac{1}{2}(x-\tilde{x})^T(x-\tilde{x}) +
  \frac{1}{2}(z-\tilde{z})^T(z-\tilde{z})$. When $K$
  satisfies~\eqref{eq:gamma_bound}, again it holds that $K \ge
  \NormInf{A^T\tilde{z}+c}$. Applying Theorem~\ref{th:grad_flow} once
  again, $\tilde{V}^{-1}(\le \rho)$ is strongly invariant
  under~\eqref{eq:flow}. Consequently, for $t \ge t_{k^*}$, we have
  $(x(t),z(t)) \in \tilde{V}^{-1}( \le
  \tilde{V}(x(t_{k^*}),z(t_{k^*}))) =
  \overline{\ball}\big((\tilde{x},\tilde{z}),
  \NormTwo{(x(t_{k^*}),z(t_{k^*})) - (\tilde{x},\tilde{z})}\big)
  \subset \overline{\ball}((\tilde{x},\tilde{z}),\epsilon)$. Since
  $\epsilon$ can be taken arbitrarily small, this implies that
  $(x(t),z(t))$ converges to the point $(\tilde{x},\tilde{z}) \in
  \pdsol$.
\end{IEEEproof}

\begin{remark}\longthmtitle{Choice of parameter $K$} \label{rem:K}
  {\rm The bound~\eqref{eq:gamma_bound} for the parameter $K$ depends
    on (i) the primal-dual solution set $\pdsol$ as well as (ii) the
    initial condition, since the result is only valid when the
    dynamics start in $V^{-1}(\le \rho)$. However, if the set~$\pdsol$
    is compact, the parameter $K$ can be chosen independently of the
    initial condition since the maximization in~\eqref{eq:gamma_bound}
    would be well defined when taken over the whole set~$\pdsol$.
    We should point out that, in
    Section~\ref{sec:discontinuous-saddle-point} we introduce a
    discontinuous version of the saddle-point dynamics which does not
    involve~$K$. \oprocend}
\end{remark}

\subsection{Discontinuous saddle-point
  dynamics}\label{sec:discontinuous-saddle-point}

Here, we propose an alternative dynamics to~\eqref{eq:flow} that does
not rely on knowledge of the parameter $K$ and also converges to the
solutions of~\eqref{eq:standard_form}-\eqref{eq:dual}. We begin by
defining the \emph{nominal flow function} $\nomflow : \reals^n_{\ge 0}
\times \reals^m \rightarrow \reals^n$ by
\begin{align*}
  \nomflow(x,z) := - c - A^T(z+Ax-b).
\end{align*}
This definition is motivated by the fact that, for $(x,z) \in
\reals_{>0}^n\times\reals^m$, the set $\partial_x L^K(x,z)$ is the
singleton $\{-\nomflow (x,z)\}$.  The \emph{discontinuous saddle-point
  dynamics} is, for $i \in \until{n}$,
\begin{subequations}\label{eq:flow2}
  \begin{align}
    \dot{x}_i &=
    \begin{cases}
      \nomflow_i(x,z), & \text{if } x_i > 0,
      \\
      \max\{0,\nomflow_i(x,z)\}, & \text{if } x_i = 0,
    \end{cases}
    \label{eq:flow2_x}
    \\
    \dot{z} &= Ax-b.
    \label{eq:flow2_z}
  \end{align}
\end{subequations}
When convenient, we use the notation $f_{\text{dis}} : \reals^{n}_{\ge
  0} \times \reals^m \rightarrow \reals^n \times \reals^m$ to refer to
the discontinuous dynamics~\eqref{eq:flow2}. Note that the
discontinuous function that defines the dynamics~\eqref{eq:flow2_x} is
simply the positive projection operator, i.e., when $x_i = 0$, it
corresponds to the projection of $\nomflow_i(x,z)$ onto $\reals_{\ge
  0}$. We understand the solutions of~\eqref{eq:flow2} in the Filippov
sense. We begin our analysis by establishing a relationship between
the Filippov set-valued map of $f_{\text{dis}}$ and the saddle-point
dynamics $\sdl$ which allows us to relate the trajectories
of~\eqref{eq:flow2} and~\eqref{eq:flow}.

\begin{proposition}\longthmtitle{Trajectories of the discontinuous
    saddle-point dynamics are trajectories of the saddle-point
    dynamics}\label{prop:filippov_subset}
  Let $\rho > 0$ and $(x_*,z_*) \in \pdsol$ be given and the function
  $V$ be defined as in Theorem~\ref{th:grad_flow}. Then, for any
  \begin{align*}
    \infty > K \ge K_1 := \max_{(x,z) \in V^{-1}(\le \rho)}
    \NormInf{\nomflow(x,z)},
  \end{align*}
  the inclusion $\mathcal{F}[f_{\text{dis}}](x,z) \subseteq \sdl(x,z)$
  holds for every $(x,z) \in V^{-1}(\le \rho)$. Thus, the trajectories
  of~\eqref{eq:flow2} starting in $V^{-1}(\le \rho)$ are also
  trajectories of~\eqref{eq:flow}.
\end{proposition}
\begin{IEEEproof}
  The projection onto the $i^{\text{th}}$ component of the Filippov
  set-valued map $\mathcal{F}[f_{\text{dis}}]$ is
  \begin{align*}
    \proj_i(\mathcal{F}[f_{\text{dis}}](x,z)) = 
    \begin{cases}
      \{\nomflow_i(x,z)\}, & \text{if $i \in \{1,\dots,n\}$ and $x_i
        > 0$},
      \\
        [\nomflow_i(x,z),\max\{0,\nomflow_i(x,z)\}], & \text{if $i \in
          \{1,\dots,n\}$ and $x_i = 0$},
        \\
        \{(Ax-b)_i\}, & \text{if $i \in \{n+1,\dots,n+m\}$.}
    \end{cases}
  \end{align*}
  As a consequence, for any $i \in \{n+1,\dots,n+m\}$, we have
  \begin{align*}
    \proj_i(\sdl(x,z)) &= (Ax-b)_i =
    \proj_i(\mathcal{F}[f_{\text{dis}}](x,z)),
  \end{align*}
  and, for any $i \in \{1,\dots,n\}$ such that $x_i > 0$, we have
  \begin{align*}
    \proj_i(\sdl(x,z)) &= (-c-A^T(Ax-b+z))_i = \{ \nomflow_i(x,z) \} =
    \proj_i(\mathcal{F}[f_{\text{dis}}](x,z)).
  \end{align*}
  Thus, let us consider the case when $x_i = 0$ for some $i \in
  \{1,\dots,n\}$. In this case, note that
  \begin{align*}
    \proj_i(\mathcal{F}[f_{\text{dis}}](x,z)) & =
    [\nomflow_i(x,z),\max\{0,\nomflow_i(x,z)\}] \subseteq
    [\nomflow_i(x,z),\nomflow_i(x,z) + |\nomflow_i(x,z)|],
    \\
    \proj_i(\sdl(x,z)) & = [\nomflow_i(x,z), \nomflow_i(x,z) + K].
  \end{align*}
  The choice $K \ge |\nomflow_i(x,z)|$ for each $i \in \{1,\dots,n\}$
  makes $\mathcal{F}[f_{\text{dis}}](x,z) \subseteq \sdl(x,z)$.  More
  generally, since $V^{-1}(\rho)$ is compact and $\nomflow$ is
  continuous, the choice
  \begin{align*}
    \infty > K \ge \max_{(x,z) \in V^{-1}(\rho)}
    \NormInf{\nomflow(x,z)},
  \end{align*}
  guarantees $\mathcal{F}[f_{\text{dis}}](x,z) \subseteq \sdl(x,z)$
  for all $(x,z) \in V^{-1}(\rho)$. By Theorem~\ref{th:grad_flow},
  we know that~$V$ is non-increasing along~\eqref{eq:flow}, implying
  that $V^{-1}(\le \rho)$ is strongly invariant with respect
  to~\eqref{eq:flow}, and hence~\eqref{eq:flow2} too. Therefore, any
  trajectory of~\eqref{eq:flow2} starting in $V^{-1}(\le \rho)$ is a
  trajectory of~\eqref{eq:flow}.
\end{IEEEproof}

\begin{figure}[hbt!]
  \centering
  \includegraphics[trim=10 50 10 30,width=.5\linewidth]{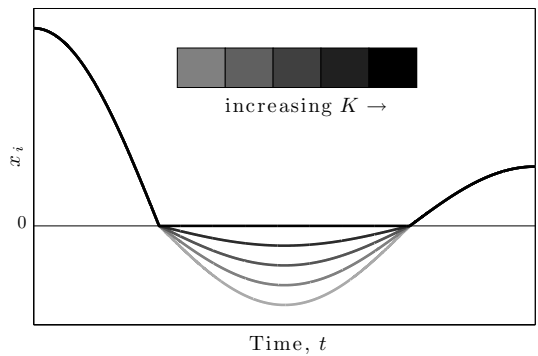}
  \caption{Illustration of the effect that increasing $K$ has
    on~\eqref{eq:flow}. For a fixed initial condition, the trajectory
    of~\eqref{eq:flow} has increasingly smaller ``incursions'' into
    the region where $x_i<0$ as $K$ increases, until a finite value is
    reached where the corresponding trajectory of~\eqref{eq:flow2} is
    also a trajectory of~\eqref{eq:flow}.}\label{fig:gamma}
\end{figure}

Note that the inclusion in Proposition~\ref{prop:filippov_subset} may
be strict and that the set of trajectories of~\eqref{eq:flow} is, in
general, richer than the set of trajectories
of~\eqref{eq:flow2}. Figure~\ref{fig:gamma} illustrates the effect
that increasing $K$ has on~\eqref{eq:flow}. From a given initial
condition, at some point the value of $K$ is large enough,
cf. Proposition~\ref{prop:filippov_subset}, to make the trajectories
of~\eqref{eq:flow2} (which never leave $\reals^n_{\ge 0} \times
\reals^m$) also be a trajectory of~\eqref{eq:flow}.

Building on Proposition~\ref{prop:filippov_subset}, the next result
characterizes the asymptotic convergence of~\eqref{eq:flow2}.

\begin{corollary}\longthmtitle{Asymptotic convergence of the
    discontinuous saddle-point dynamics}\label{cor:grad_flow_modified}
  The trajectories of~\eqref{eq:flow2} starting in $\reals_{\ge 0}^n
  \times \reals^m$ converge asymptotically to a point in $\pdsol$.
\end{corollary}
\begin{IEEEproof}
  Let $V$ be defined as in Theorem~\ref{th:grad_flow}. Given any
  initial condition $(x_0,z_0) \in \reals^n \times \reals^m$, let $t
  \mapsto (x(t),z(t))$ be a trajectory of~\eqref{eq:flow2} starting
  from $(x_0,z_0)$ and let $\rho = V(x_0,z_0)$. Note that $t \mapsto
  (x(t),z(t))$ does not depend on $K$ because~\eqref{eq:flow2} does
  not depend on $K$.  Proposition~\ref{prop:filippov_subset}
  establishes that $t \mapsto (x(t),z(t))$ is also a trajectory
  of~\eqref{eq:flow} for $K \ge K_1$. Imposing the additional
  condition that
  \begin{align*}
    \infty > K \ge \max \bigg\{ K_1, \max_{(x_*,z_*) \in (\pdsol) \cap
      V^{-1}(\le \rho)}\NormInf{A^Tz_* + c} \bigg \},
  \end{align*}
  Corollary~\ref{cor:point} implies that the trajectories
  of~\eqref{eq:flow} (and thus $t \mapsto (x(t),z(t)$) converge
  asymptotically to a point in $\pdsol$.
\end{IEEEproof}

One can also approach the convergence analysis of~\eqref{eq:flow2}
from a switched systems perspective, which would require checking that
certain regularity conditions hold for the switching behavior of the
system. We have been able to circumvent this complexity by relying on
the powerful stability tools available for set-valued dynamics to
analyze~\eqref{eq:flow} and by relating its solutions with those
of~\eqref{eq:flow2}.
Moreover, the interpretation of the trajectories of~\eqref{eq:flow2}
in the Filippov sense is instrumental for our analysis in
Section~\ref{sec:robustness} where we study the robustness against
disturbances using powerful Lyapunov-like tools for differential
inclusions.

\begin{remark}\longthmtitle{Comparison to existing dynamics for linear
    programming} 
  {\rm Though a central motivation for the development of our linear
    programming algorithm is the establishment of various robustness
    properties which we study next, the dynamics~\eqref{eq:flow2} and
    associated convergence results of this section are both novel and
    have distinct contributions. The work~\cite{KA-LH-HU:58} builds on
    the saddle-point dynamics of a smooth Lagrangian function to
    introduce an algorithm for linear programming.  Instead of exact
    penalty functions, this approach uses projections to keep the
    evolution within the feasible set, resulting in a discontinuous
    dynamics in both the primal and dual variables.  The
    work~\cite{DF-FP:10} employs a similar approach to deal with
    non-strictly convex programs under inequality constraints, where
    projection is used instead employed to keep nonnegative the value
    of the dual variables.  These works establish convergence in the
    primal variables (\cite{KA-LH-HU:58} under the assumption that the
    solution of the linear program is unique,~\cite{DF-FP:10} under the
    assumption that Slater's condition is satisfied) to a solution of
    the linear program. In both cases, the dual variables converge to
    some unknown point which might not be a solution to the dual
    problem.  This is to be contrasted with the convergence properties
    of the dynamics~\eqref{eq:flow2} stated in
    Corollary~\ref{cor:grad_flow_modified} which only require the
    linear program to be feasible with finite optimal value. }
  \oprocend

\end{remark}

\subsection{Distributed implementation}\label{sec:DLP_dist}

An important advantage of the dynamics~\eqref{eq:flow2} over other
linear programming methods is that it is well-suited for distributed
implementation. To make this statement precise, consider a scenario
where each component of $x \in \reals^n$ corresponds to an independent
decision maker or agent and the interconnection between the agents is
modeled by an undirected graph $\graph = (\vertices,\edges)$. To see
under what conditions the dynamics~\eqref{eq:flow2} can be implemented
by this multi-agent system, let us express it component-wise. First,
the nominal flow function in~\eqref{eq:flow2_x} for agent $i \in
\until{n}$ is,
\begin{align*}
  \nomflow_i(x,z) = - c_i - \sum_{\ell = 1}^m a_{\ell,i} \Big[z_{\ell}
  + \sum_{k=1}^n a_{\ell,k}x_k - b_\ell\Big] = - c_i - \hspace{-4.7mm}
  \sum_{\{\ell \; : \; a_{\ell,i} \not= 0\}} \hspace{-4.7mm}
  a_{\ell,i} \Big[z_{\ell} + \hspace{-5.1mm} \sum_{\{k \; : \;
    a_{\ell,k} \not= 0 \}} \hspace{-5.1mm} a_{\ell,k}x_k -
  b_\ell\Big], \nonumber
\end{align*}
and the dynamics~\eqref{eq:flow2_z} for each $\ell \in \{1,\dots,m\}$
is
\begin{align}\label{eq:flow_z_ell}
  \dot{z}_{\ell} = \sum_{\{ i \; : \; a_{\ell,i} \not= 0\}}
  a_{\ell,i}x_i - b_{\ell}.
\end{align} 
According to these expressions, in order for agent $i \in \until{n}$
to be able to implement its corresponding dynamics
in~\eqref{eq:flow2_x}, it also needs access to certain components of
$z$ (specifically, those components $z_\ell$ for which $a_{\ell,i}
\neq 0$), and therefore needs to implement their corresponding
dynamics~\eqref{eq:flow_z_ell}. We say that the
dynamics~\eqref{eq:flow2} is \emph{distributed over $\graph$} when the
following holds
\begin{enumerate}
\item[(D1)] for each $i \in \vertices$, agent $i$ knows
  \begin{enumerate}
  \item $c_i \in \reals$,
  \item every $b_{\ell} \in \reals$ for which $a_{\ell,i} \not=
    0$,
  \item the non-zero elements of every row of $A$ for which the
    $i^{\text{th}}$ component, $a_{\ell,i}$, is non-zero,
  \end{enumerate}
\item[(D2)] agent $i \in \vertices$ has control over the variable $x_i
  \in \reals$,
\item[(D3)] $\graph$ is connected with respect to $A$, and
\item[(D4)] agents have access to the variables controlled by
  neighboring agents.
\end{enumerate}
Note that (D3) guarantees that the agents that
implement~\eqref{eq:flow_z_ell} for a particular $\ell \in \until{m}$
are neighbors in $\graph$. 

\begin{remark}\longthmtitle{Scalability of the nominal saddle-point
    dynamics} {\rm A different approach to
    solve~\eqref{eq:standard_form} is the following: reformulate the
    optimization problem as the constrained minimization of a sum of
    convex functions all of the form $\frac{1}{n} c^T x$ and use the
    algorithms developed in, for
    instance,~\cite{AN-AO-PAP:10,MZ-SM:12,MB-GN-FB-FA:12,GN-FB:11,JW-NE:11},
    for distributed convex optimization. However, in this case, this
    approach would lead to agents storing and communicating with
    neighbors estimates of the entire solution vector in $\reals^n$,
    and hence would not scale well with the number of agents of the
    network. In contrast, to execute the discontinuous saddle-point
    dynamics, agents only need to store the component of the solution
    vector that they control and communicate it with
    neighbors. Therefore, the dynamics scales well with respect to the
    number of agents in the network.
  }\oprocend
\end{remark}

\myclearpage
\section{Robustness against disturbances}\label{sec:robustness}

Here we explore the robustness properties of the discontinuous
saddle-point dynamics~\eqref{eq:flow2} against disturbances. Such
disturbances may correspond to noise, unmodeled dynamics, or incorrect
agent knowledge of the data defining the linear program. Note that the
global asymptotic stability of $\pdsol$ under~\eqref{eq:flow2}
characterized in Section~\ref{sec:DLP} naturally provides a robustness
guarantee on this dynamics: when $\pdsol$ is compact, sufficiently
small perturbations do not destroy the global asymptotic stability of
the equilibria, cf.~\cite{CC-ART-RG:08}. Our objective here is to go
beyond this qualitative statement to obtain a more precise,
quantitative description of robustness.  To this end, we consider the
notions of input-to-state stability (ISS) and integral-input-to-state
stability (iISS).  In Section~\ref{sec:not-ISS} we show that, when the
disturbances correspond to uncertainty in the problem data, no
dynamics for linear programming can be ISS. This motivates us to
explore the weaker notion of iISS. In Section~\ref{se:saddle-are-iISS}
we show that~\eqref{eq:flow2} with additive disturbances is~iISS.

\begin{remark}\longthmtitle{Robust dynamics versus robust
    optimization}
  {\rm We make a note of the distinction between the notion of
    algorithm robustness, which is what we study here, and the term
    robust (or worst-case) optimization, see
    e.g.,~\cite{DB-DBB-CC:11}. 
    The latter refers to a type of problem formulation in which some
    notion of variability (which models uncertainty) is explicitly
    included in the problem statement. Mathematically,
    \begin{align*}
      \min \; \; c^Tx \quad \text{s.t.} \; \; f(x,\omega) \le 0, \;
      \forall \omega \in \Omega,
    \end{align*}
    where $\omega$ is an uncertain parameter. Building on the
    observation that one only has to consider the worst-case values of
    $\omega$, one can equivalently cast the optimization problem with
    constraints that only depend on $x$, albeit at the cost of a loss
    of structure in the formulation.  Another point of connection with
    the present work is the body of research on stochastic
    approximation in discrete optimization, where the optimization
    parameters are corrupted by disturbances, see
    e.g.~\cite{HJK-GGY:03}.
    \oprocend }
\end{remark}

Without explicitly stating it from here on, we make the following
assumption along the section:
\begin{enumerate}
\item[{\bf (A)}] The solution sets to~\eqref{eq:standard_form}
  and~\eqref{eq:dual} are compact (i.e., $\pdsol$ is compact).
\end{enumerate}
The justification for this assumption is twofold. On the technical
side, our study of the iISS properties of~\eqref{eq:flow_robust} in
Section~\ref{se:saddle-are-iISS} builds on a Converse Lyapunov
Theorem~\cite{CC-ART-RG:08} which requires the equilibrium set to be
compact (the question of whether the Converse Lyapunov Theorem holds
when the equilibrium set is not compact and the dynamics is
discontinuous is an open problem).
On the practical side, one can add box-type constraints
to~\eqref{eq:standard_form}, ensuring that (A) holds.

We now formalize the disturbance model considered in this
section. Let $w = (w_x,w_z): \realsnonnegative \rightarrow \reals^n
\times \reals^m$ be locally essentially bounded and enter the dynamics
as follows,
\begin{subequations}\label{eq:flow_robust}
  \begin{align}
    \dot{x}_i &=
    \begin{cases}
      \nomflow_i(x,z) + (w_x)_i, & \text{if } x_i > 0,
      \\
      \max\{0,\nomflow_i(x,z) + (w_x)_i\}, & \text{if } x_i = 0,
    \end{cases} \hspace{6mm} \forall i \in
    \{1,\dots,n\},
    \label{eq:flow_robust_x}
    \\
    \dot{z} &= Ax-b + w_z.
    \label{eq:flow_robust_z}
  \end{align}
\end{subequations}
For notational purposes, we use $\sdlw: \reals^{2(n+m)} \rightarrow
\reals^{n+m}$ to denote~\eqref{eq:flow_robust}. We exploit the fact
that $\nomflow$ is affine to state that the additive disturbance $w$
captures unmodeled dynamics, measurement and computation noise, and
any error in an agent's knowledge of the problem data ($A,b$ and
$c$). For example, if agent $i \in \until{n}$ uses an estimate
$\hat{c}_i$ of $c_i$ when computing its dynamics, this can be modeled
in~\eqref{eq:flow_robust} by considering $(w_{x}(t))_i = c_i -
\hat{c}_i$. To make precise the correspondence between the disturbance
$w$ and uncertainties in the problem data, we provide the following
convergence result when the disturbance is constant.


\begin{corollary}\longthmtitle{Convergence under constant
    disturbances}\label{cor:pert_prog} 
  For constant $\overline{w} = (\overline{w}_x,\overline{w}_z) \in
  \reals^n \times \reals^m$, consider the \emph{perturbed linear
    program},
  \begin{subequations}\label{eq:pert_prog}
    \begin{alignat}{2} 
      &\min && \quad (c - \overline{w}_x - A^T \overline{w}_z)^Tx
      \\ &\hspace{1.5mm}\text{\rm s.t.} && \quad Ax = b - \overline{w}_z,
      \quad x \ge 0,
    \end{alignat}
  \end{subequations}
  and, with a slight abuse in notation, let $\primalsol(\overline{w})
  \times \dualsol(\overline{w})$ be its primal-dual solution
  set. Suppose that $\primalsol(\overline{w}) \times
  \dualsol(\overline{w})$ is nonempty. Then each trajectory
  of~\eqref{eq:flow_robust} starting in $\reals_{\ge 0}^n \times
  \reals^m$ with constant disturbance $w(t) = \overline{w} =
  (\overline{w}_x,\overline{w}_z)$ converges asymptotically to a point
  in $\primalsol(\overline{w}) \times \dualsol(\overline{w})$.
\end{corollary}
\begin{IEEEproof}
  Note that~\eqref{eq:flow_robust} with disturbance $\overline{w}$
  corresponds to the undisturbed dynamics~\eqref{eq:flow2} for the
  perturbed problem~\eqref{eq:pert_prog}.
  Since $\primalsol(\overline{w}) \times \dualsol(\overline{w}) \neq
  \emptyset$, Corollary~\ref{cor:grad_flow_modified} implies the
  result.
\end{IEEEproof}


\subsection{No dynamics for linear programming is input-to-state
  stable}\label{sec:not-ISS}

The notion of input-to-state stability (ISS) is a natural starting
point to study the robustness of dynamical systems against
disturbances.  Informally, if a dynamics is ISS, then bounded
disturbances give rise to bounded deviations from the equilibrium set.
Here we show that any dynamics that (i) solve any feasible linear
program and (ii) where uncertainties in the problem data ($A,b$, and
$c$) enter as disturbances is not input-to-state stable (ISS). Our
analysis relies on the properties of the solution set of a linear
program.
To make our discussion precise, we begin by recalling the definition
of input-to-state stability.


\begin{definition}\longthmtitle{Input-to-state
    stability~\cite{EDS:89-tac}}
  The dynamics~\eqref{eq:flow_robust} is ISS with respect to $\pdsol$
  if there exist $\beta \in \mathcal{KL}$ and $\gamma \in \mathcal{K}$
  such that, for any trajectory $t\mapsto (x(t),z(t))$
  of~\eqref{eq:flow_robust}, one has
  \begin{align*}
    \NormTwo{(x(t),z(t))}_{\pdsol} \le
    \beta(\NormTwo{(x(0),z(0)}_{\pdsol},t) +
    \gamma(\NormInf{w}),
  \end{align*}
  for all $t \ge 0$. Here, $\NormInf{w} := \esssup_{s \ge 0}
  \NormTwo{w(s)}$ is the essential supremum of $w(t)$.
\end{definition}

Our method to show that no dynamics is ISS is constructive. We find a
constant disturbance such that the primal-dual solution set to some
perturbed linear program is unbounded. Since any point in this
unbounded solution set is a stable equilibrium by assumption, this
precludes the possibility of the dynamics from being ISS.  This
argument is made precise next.

\begin{theorem}\longthmtitle{No dynamics for linear programming is
    ISS}\label{th:no-ISS} 
  Consider the generic dynamics
  \begin{align}
    (\dot{x},\dot{z}) = \Phi(x,z,v) \label{eq:generic_dyn}
  \end{align}
  with disturbance $t \mapsto v(t)$. Assume uncertainties in the
  problem data are modeled by $v$. That is, there
  exists a surjective function $g = (g_1,g_2) : \reals^{n +m}
  \rightarrow \reals^n \times \reals^m$ with $g(0) = (0,0)$ such that,
  for $\bar{v} \in \reals^{n +m}$, the primal-dual solution set
  $\mathcal{X}(\bar{v}) \times \mathcal{Z}(\bar{v})$ of the linear
  program
  \begin{subequations}\label{eq:no_ISS_prog}
    \begin{alignat}{2}
      &\min && \quad (c+g_1(\bar{v}))^Tx
      \\
      &\hspace{1.5mm}\text{\rm s.t.} && \quad Ax = b + g_2(\bar{v}), \quad
      x \ge 0.
    \end{alignat}
  \end{subequations}
  is the stable equilibrium set of $(\dot{x},\dot{z}) =
  \Phi(x,z,\bar{v})$ whenever 
  $\mathcal{X}(\bar{v}) \times \mathcal{Z}(\bar{v}) \not= \emptyset$.
  Then, the dynamics~\eqref{eq:generic_dyn} is not ISS with respect to
  $\pdsol$.
\end{theorem}
\begin{IEEEproof}
  We divide the proof in two cases depending on whether $\{ Ax = b, x
  \ge 0 \}$ is (i) unbounded or (ii) bounded. In both cases, we design
  a constant disturbance $v(t) = \bar{v}$ such that the equilibria
  of~\eqref{eq:generic_dyn} contains points arbitrarily far away from
  $\pdsol$.  This would imply that the dynamics is not ISS. Consider
  case (i). Since $\{ Ax = b, x \ge 0 \}$ is unbounded, convex, and
  polyhedral, there exists a point $\hat{x} \in \reals^n$ and
  direction $\nu_x \in \reals^n \setminus \{0\}$ such that $\hat{x} +
  \lambda \nu_x \in \bd(\{ Ax = b, x \ge 0 \})$ for all $\lambda \ge
  0$. Here $\bd(\cdot)$ refers to the boundary of the set. Let $\eta
  \in \reals^n$ be such that $\eta^T \nu_x = 0$ and $\hat{x} +
  \epsilon \eta \notin \{ Ax = b, x \ge 0 \}$ for any $\epsilon > 0$
  (geometrically, $\eta$ is normal to and points out of $\{ Ax = b, x
  \ge 0 \}$ at $\hat{x}$).  Now that these quantities have been
  defined, consider the following linear program,
  \begin{align} \label{eq:iss_contra}
    \min \quad \eta^Tx \quad \text{s.t.} \quad Ax = b, \quad x \ge 0.
  \end{align}
  Because $g$ is surjective, there exists $\bar{v}$ such that
  $g(\bar{v}) = (-c+\eta,0)$. In this case, the
  program~\eqref{eq:iss_contra} is exactly the
  program~\eqref{eq:no_ISS_prog}, with primal-dual solution set
  $\mathcal{X}(\bar{v}) \times \mathcal{Z}(\bar{v})$. We show next
  that $\hat{x}$ is a solution to~\eqref{eq:iss_contra} and thus in
  $\mathcal{X}(\bar{v})$. Clearly, $\hat{x}$ satisfies the constraints
  of~\eqref{eq:iss_contra}.  Since $\eta^T \nu_x = 0$ and points
  outward of $\{ Ax = b, x \ge 0 \}$, it must be that $\eta^T(\hat{x}
  - x) \le 0$ for any $x \in \{ Ax = b, x \ge 0 \}$, which implies
  that $\eta^T\hat{x} \le \eta^T x$. Thus, $\hat{x}$ is a solution
  to~\eqref{eq:iss_contra}.  Moreover, $\hat{x} + \lambda \nu_x$ is
  also a solution to~\eqref{eq:iss_contra} for any $\lambda \ge 0$
  since (i) $\eta^T (\hat{x} + \lambda \nu_x) = \eta^T\hat{x}$ and
  (ii) $\hat{x} + \lambda \nu_x \in \{ Ax = b, x \ge 0 \}$. That is,
  $\mathcal{X}(\bar{v})$ is unbounded. Therefore, there is a point
  $(x_0,z_0) \in \mathcal{X}(\bar{v}) \times \mathcal{Z}(\bar{v})$,
  which is also an equilibrium of~\eqref{eq:generic_dyn} by
  assumption, that is arbitrarily far from the set $\pdsol$. Clearly,
  $t \mapsto (x(t),z(t)) = (x_0,z_0)$ is an equilibrium trajectory
  of~\eqref{eq:generic_dyn} starting from $(x_0,z_0)$ when $v(t) =
  \bar{v}$. The fact that $(x_0,z_0)$ can be made arbitrarily far from
  $\pdsol$ precludes the possibility of the dynamics from being ISS.

  Next, we deal with case (ii), when $\{ Ax = b, x \ge 0 \}$ is
  bounded.
  Consider the linear program
  \begin{align*}
    \max \quad -b^Tz \quad \text{s.t.} \quad A^Tz \ge 0.
  \end{align*}
  Since $\{Ax = b,x \ge 0 \}$ is bounded, Lemma~\ref{lem:prop-feas}
  implies that $\{ A^Tz \ge 0 \}$ is unbounded.  Using an
  analogous approach as in case (i), one can find $\eta \in
  \reals^m$ such that the set of solutions to
  \begin{align}
    \max \quad \eta^Tz \quad \text{s.t.} \quad A^Tz
    \ge 0, \label{eq:unbounded-dual}
  \end{align}
  is unbounded.  Because $g$ is surjective, there exists $\bar{v}$
  such that $g(\bar{v}) = (-c,-b-\eta)$. In this case, the
  program~\eqref{eq:unbounded-dual} is the dual
  to~\eqref{eq:no_ISS_prog}, with primal-dual solution set
  $\mathcal{X}(\bar{v}) \times \mathcal{Z}(\bar{v})$. Since
  $\mathcal{Z}(\bar{v})$ is unbounded, one can find equilibrium
  trajectories of~\eqref{eq:generic_dyn} under the disturbance $v(t) =
  \bar{v}$ that are arbitrarily far away from $\pdsol$, which
  contradicts ISS.
\end{IEEEproof}

Note that, in particular, the perturbed problem~\eqref{eq:pert_prog}
and~\eqref{eq:no_ISS_prog} coincide when
\begin{align*}
  g(\overline{w}) = g(\overline{w}_x,\overline{w}_z) =
  (-\overline{w}_x - A^T\overline{w}_z, -\overline{w}_z).
\end{align*}
Thus, by Theorem~\ref{th:no-ISS}, the discontinuous saddle-point
dynamics~\eqref{eq:flow_robust} is not ISS.  Nevertheless, one can
establish an ISS-like result for this dynamics under small enough and
constant disturbances. We state this result next, where we also
provide a quantifiable upper bound on the disturbances in terms of the
solution set of some perturbed linear program.

\begin{proposition}\longthmtitle{ISS of discontinuous saddle-point
    dynamics under small constant
    disturbances}\label{prop:small_constant} 
  Suppose there exists $\delta > 0$ such that the primal-dual solution
  set $\primalsol(\overline{w}) \times \dualsol(\overline{w})$ of the
  perturbed problem~\eqref{eq:pert_prog} is nonempty for $\overline{w}
  \in \overline{\ball}(0,\delta)$ and $\cup_{\overline{w} \in
    \overline{\ball}(0,\delta)} \primalsol(\overline{w}) \times
  \dualsol(\overline{w})$ is compact. Then there exists a continuous,
  zero-at-zero, and increasing function $\gamma: [0,\delta]
  \rightarrow \reals_{\ge 0}$ such that, for all trajectories
  $t\mapsto (x(t),z(t))$ of~\eqref{eq:flow_robust} with constant
  disturbance $\overline{w} \in \overline{\ball}(0,\delta)$, it holds
  that
  \begin{align*}
    \lim_{t \rightarrow \infty} \NormTwo{(x(t),z(t))}_{\pdsol} \le
    \gamma(\NormTwo{\overline{w}}).
  \end{align*}
\end{proposition}
\begin{IEEEproof}
  Let $\gamma: [0,\delta] \rightarrow \reals_{\ge 0}$ be given by
  \begin{align*}
    \gamma (r) := \max \bigg \{ \NormTwo{(x,z)}_{\pdsol} : (x,z) \in
    \bigcup_{\overline{w} \in \overline{\ball}(0,r)}
    \primalsol(\overline{w}) \times \dualsol(\overline{w}) \bigg \} .
  \end{align*}
  By hypotheses, $\gamma$ is well-defined.  Note also that $\gamma$ is
  increasing and satisfies $\gamma(0) = 0$. Next, we show that
  $\gamma$ is continuous. By assumption, $\primalsol(\overline{w})
  \times \dualsol(\overline{w})$ is nonempty and bounded for every
  $\overline{w} \in \overline{\ball}(0,\delta)$. Moreover, it is clear
  that $\primalsol(\overline{w}) \times \dualsol(\overline{w})$ is
  closed for every $\overline{w} \in \overline{\ball}(0,\delta)$ since
  we are considering linear programs in standard form. Thus,
  $\primalsol(\overline{w}) \times \dualsol(\overline{w})$ is nonempty
  and compact for every $\overline{w} \in \overline{\ball}(0,\delta)$.
  By~\cite[Corollary 11]{RJBW:85}, these two conditions are sufficient
  for the set-valued map $\overline{w} \mapsto
  \primalsol(\overline{w}) \times \dualsol(\overline{w})$ to be
  continuous on $\overline{\ball}(0,\delta)$. Since $r \mapsto
  \overline{\ball}(0,r)$ is also continuous,~\cite[Proposition 1,
  pp. 41]{JPA-AC:84} ensures that the following set-valued composition
  map
  \begin{align*}
    r \mapsto \bigcup_{\overline{w} \in \overline{\ball}(0,r)}
    \primalsol(\overline{w}) \times \dualsol(\overline{w})
  \end{align*}
  is continuous (with compact values, by
  assumption). Therefore,~\cite[Theorem 6, pp. 53]{JPA-AC:84}
  guarantees then that $\gamma$ is continuous on
  $\overline{\ball}(0,\delta)$. Finally, to establish the bound on the
  trajectories, recall from Corollary~\ref{cor:pert_prog} that each
  trajectory $t \mapsto (x(t),z(t))$ of~\eqref{eq:flow_robust} with
  constant disturbance $\overline{w} \in \overline{\ball}(0,\delta)$
  converges asymptotically to a point in $\primalsol(\overline{w})
  \times \dualsol(\overline{w})$.  The distance between $\pdsol$ and
  the point in $\primalsol(\overline{w}) \times
  \dualsol(\overline{w})$ to which the trajectory converges is upper
  bounded by
  \begin{align*}
    \lim_{t \rightarrow \infty} \NormTwo{(x(t),z(t))}_{\pdsol} &\le
    \max \{ \NormTwo{(x,z)}_{\pdsol} : (x,z) \in
    \primalsol(\overline{w}) \times \dualsol(\overline{w}) \} \le
    \gamma(\NormTwo{\overline{w}}) ,
  \end{align*}
  which concludes the proof.
\end{IEEEproof}

\subsection{Discontinuous saddle-point dynamics is integral
  input-to-state stable}\label{se:saddle-are-iISS}

Here we establish that the dynamics~\eqref{eq:flow_robust} possess a
notion of robustness weaker than ISS, namely, integral input-to-state
stability (iISS). Informally, iISS guarantees that disturbances with
small energy give rise to small deviations from the equilibria. This
is stated formally next.

\begin{definition}\longthmtitle{Integral input-to-state
    stability~\cite{DA-EDS-YW:00}}
  The dynamics~\eqref{eq:flow_robust} is iISS with respect to the set
  $\pdsol$ if there exist functions $\alpha \in \mathcal{K}_{\infty},
  \beta \in \mathcal{KL},$ and $\gamma \in \mathcal{K}$ such that, for
  any trajectory $t\mapsto (x(t),z(t))$ of~\eqref{eq:flow_robust} and
  all $t \ge 0$, one has
  \begin{align}\label{eq:iISS_characterization}
    \alpha(\NormTwo{(x(t),z(t))}_{\pdsol}) \le
    \beta(\NormTwo{(x(0),z(0)}_{\pdsol},t) + \int_0^t
    \gamma(\NormTwo{w(s)})ds.
  \end{align}
\end{definition}

Our ensuing discussion is based on a suitable adaptation of the
exposition in~\cite{DA-EDS-YW:00} to the setup of asymptotically
stable sets for discontinuous dynamics. A useful tool for establishing
iISS is the notion of iISS Lyapunov function, whose definition we
review next.

\begin{definition}\longthmtitle{iISS Lyapunov function}
  A differentiable function $V:\reals^{n+m} \rightarrow \reals_{\ge
    0}$ is an iISS Lyapunov function with respect to the set
  $\pdsol$ for dynamics~\eqref{eq:flow_robust} if
  there exist functions $\alpha_1,\alpha_2 \in \mathcal{K}_{\infty},
  \sigma \in \mathcal{K},$ and a continuous positive definite function
  $\alpha_3$ such that
  \begin{subequations}
    \begin{gather}
      \alpha_1(\NormTwo{(x,z)}_{\pdsol}) \le V(x,z)
      \le \alpha_2(\NormTwo{(x,z)}_{\primalsol \times
        \dualsol}), \label{eq:proper}
      \\
      a \le - \alpha_3(\NormTwo{(x,z)}_{\pdsol}) +
      \sigma(\NormTwo{w}), \label{eq:dissipate}
    \end{gather}    
  \end{subequations}
  for all $a \in \mathcal{L}_{\mathcal{F}[\sdlw]}V(x,z)$ and $x \in
  \reals^n, z \in \reals^m,w \in \reals^{n+m}$.
\end{definition}

Note that, since the set $\pdsol$ is compact (cf. Assumption
(A)),~\eqref{eq:proper} is equivalent to $V$ being proper with respect
to~$\pdsol$.  The existence of an iISS Lyapunov function is critical
in establishing iISS, as the following result states.

\begin{theorem}\longthmtitle{iISS Lyapunov function implies
    iISS}\label{th:iISS_lyap}
  If there exists an iISS Lyapunov function with respect to $\pdsol$
  for~\eqref{eq:flow_robust}, then the dynamics is iISS with respect
  to $\pdsol$.
\end{theorem}

This result is stated in~\cite[Theorem 1]{DA-EDS-YW:00} for
the case of differential equations with locally Lipschitz right-hand
side and asymptotically stable origin, but its extension to
discontinuous dynamics and asymptotically stable sets, as considered
here, is straightforward. We rely on Theorem~\ref{th:iISS_lyap} to
establish that the discontinuous saddle-point
dynamics~\eqref{eq:flow_robust} is iISS. Interestingly, the function
$V$ employed to characterize the convergence properties of the
unperturbed dynamics in Section~\ref{sec:DLP} is not an iISS Lyapunov
function (in fact, our proof of Theorem~\ref{th:grad_flow} relies on
the set-valued LaSalle Invariance Principle because, essentially, the
Lie derivative of $V$ is not negative definite). Nevertheless, in the
proof of the next result, we build on the properties of this function
with respect to the dynamics to identify a suitable iISS Lyapunov
function for~\eqref{eq:flow_robust}.

\begin{theorem}\longthmtitle{iISS of saddle-point
    dynamics}\label{th:saddle_iISS}
  The dynamics~\eqref{eq:flow_robust} is iISS with respect to
  $\primalsol \times \dualsol$.
\end{theorem}
\begin{IEEEproof}
  We proceed by progressively defining functions $\Vcirc$,
  $\Vcircsep$, $\Vclf$, and $\Vclfsep: \reals^n \times \reals^m
  \rightarrow \reals$. The rationale for our construction is as
  follows. Our starting point is the squared Euclidean distance from
  the primal-dual solution set, denoted $\Vcirc$.  The function
  $\Vcircsep$ is a reparameterization of~$\Vcirc$ (which remains
  radially unbounded with respect to $\pdsol$) so that state and
  disturbance appear separately in the (set-valued) Lie
  derivative. However, since $\Vcirc$ is only a LaSalle-type function,
  this implies that only the disturbance appears in the Lie derivative
  of~$\Vcircsep$. Nevertheless, via a Converse Lyapunov Theorem, we
  identify an additional function $\Vclf$ whose reparameterization
  $\Vclfsep$ has a Lie derivative where both state and disturbance
  appear. The function $\Vclfsep$, however, may not be radially
  unbounded with respect to $\pdsol$. This leads us to the
  construction of the iISS Lyapunov function as $V = \Vcircsep +
  \Vclfsep$.

  We begin by defining the differentiable function $\Vcirc$
  \begin{align*}
    \Vcirc(x,z) = \min_{(x_*,z_*) \in \pdsol} \frac{1}{2}(x-x_*)^T(x-x_*)
    + \frac{1}{2}(z-z_*)^T(z-z_*).
  \end{align*}
  Since $\pdsol$ is convex and compact, applying
  Theorem~\ref{th:danskin} one gets $\nabla \Vcirc(x,z) =
  (x-x_*(x,z),z-z_*(x,z))$, where
  \begin{align*}
    (x_*(x,z), z_*(x,z)) = \underset{(x_*,z_*) \in \primalsol \times
      \dualsol}{\argmin} \frac{1}{2}(x-x_*)^T(x-x_*) +
    \frac{1}{2}(z-z_*)^T(z-z_*).
  \end{align*}
  It follows from Theorem~\ref{th:grad_flow} and
  Proposition~\ref{prop:filippov_subset} that
  $\mathcal{L}_{\mathcal{F}[f_{\text{dis}}]}\Vcirc(x,z) \subset
  (-\infty,0]$ for all $(x,z) \in \reals^n_{\ge 0} \times
  \reals^m$. Next, similar to the approach in~\cite{DA-EDS-YW:00},
  define the function $\Vcircsep$ by
  \begin{align*}
    \Vcircsep(x,z) = \scalebox{1.2}{$\int_0^{\Vcirc(x,z)}
      \frac{dr}{1+\sqrt{2r}}$}.
  \end{align*}
  Clearly, $\Vcircsep(x,z)$ is positive definite with respect to
  $\primalsol \times \dualsol$. Also, $\Vcircsep(x,z)$ is radially
  unbounded with respect to $\pdsol$ because (i) $\Vcirc(x,z)$ is
  radially unbounded with respect to $\pdsol$ and (ii) $\lim_{y
    \rightarrow \infty} \int_0^y \frac{dr}{1+\sqrt{2r}} = \infty$. In
  addition, for any $a \in
  \mathcal{L}_{\mathcal{F}[\sdlw]}\Vcircsep(x,z)$ and $(x,z) \in
  \reals^n_{\ge 0} \times \reals^m$, one has
  \begin{align}\label{eq:W1_1}
    a &\le \frac{\sqrt{2\Vcirc(x,z)}
      \NormTwo{w}}{1+\sqrt{2\Vcirc(x,z)}} \le \NormTwo{w}.
  \end{align}
  Next, we define the function $\Vclf$. Since $\primalsol \times
  \dualsol$ is compact and globally asymptotically stable
  for~\eqref{eq:flow2} $(\dot{x},\dot{z}) = \mathcal{F}[\sdlw](x,z)$
  when $w \equiv 0$ (cf.  Corollary~\ref{cor:grad_flow_modified}) the
  Converse Lyapunov Theorem~\cite[Theorem 3.13]{CC-ART-RG:08} ensures
  the existence of a smooth function $\Vclf:\reals^{n+m} \rightarrow
  \reals_{\ge 0}$ and class $\mathcal{K}_{\infty}$ functions
  $\tilde{\alpha}_1$, $\tilde{\alpha}_2$, $\tilde{\alpha}_3$ such that
  \begin{gather*}
    \tilde{\alpha}_1(\NormTwo{(x,z)}_{\pdsol}) \le \Vclf(x,z) \le
    \tilde{\alpha}_2(\NormTwo{(x,z)}_{\pdsol}) ,
    \\
    \quad a \le -\tilde{\alpha}_3(\NormTwo{(x,z)}_{\primalsol \times
      \dualsol}) ,
  \end{gather*}
  for all $a \in \mathcal{L}_{\mathcal{F}[f_{\text{dis}}]}\Vclf(x,z)$
  and $(x,z) \in \reals^n_{\ge 0} \times \reals^m$.  Thus, when $w
  \not \equiv 0$, for $a \in
  \mathcal{L}_{\mathcal{F}[\sdlw]}\Vclf(x,z)$ and $(x,z) \in
  \reals^n_{\ge 0} \times \reals^m$, we have
  \begin{align*}
    a &\le -\tilde{\alpha}_3(\NormTwo{(x,z)}_{\primalsol \times
      \dualsol}) + \nabla \Vclf(x,z) w,
    \\
    &\le -\tilde{\alpha}_3(\NormTwo{(x,z)}_{\primalsol \times
      \dualsol}) + \NormTwo{\nabla \Vclf(x,z)} \cdot \NormTwo{w},
    \\
    &\le -\tilde{\alpha}_3(\NormTwo{(x,z)}_{\primalsol \times
      \dualsol})  + (\NormTwo{(x,z)}_{\pdsol} +
    \NormTwo{\nabla \Vclf(x,z)} ) \cdot \NormTwo{w},
    \\
    &\le -\tilde{\alpha}_3(\NormTwo{(x,z)}_{\primalsol \times
      \dualsol}) + \lambda(\NormTwo{(x,z)}_{\primalsol \times
      \dualsol}) \cdot \NormTwo{w} ,
  \end{align*}
  where $\lambda : [0,\infty) \rightarrow [0,\infty)$ is given by
  \begin{align*}
    \lambda(r) = r + \max_{\NormTwo{\eta}_{\pdsol}
      \le r} \NormTwo{\nabla \Vclf(\eta)}.
  \end{align*}
  Since $\Vclf$ is smooth, $\lambda$ is a class $\mathcal{K}$
  function. Next, define
  \begin{align*}
    \Vclfsep(x,z) = \scalebox{1.2}{$\int_0^{\Vclf(x,z)}
    \frac{dr}{1+\lambda \circ \tilde{\alpha}_1^{-1}(r)}$}.
  \end{align*}
  Without additional information about $\lambda \circ
  \tilde{\alpha}_1^{-1}$, one cannot determine if $\Vclfsep$ is
  radially unbounded with respect to $\pdsol$ or not. Nevertheless,
  $\Vclfsep$ is positive definite with respect to $\pdsol$. Then for
  any $a \in \mathcal{L}_{\mathcal{F}[\sdlw]}\Vclfsep(x,z)$ and $(x,z)
  \in \reals^n_{\ge 0} \times \reals^m$ we have,
  \begin{align}
    a &\le \frac{-\tilde{\alpha}_3(\NormTwo{(x,z)}_{\primalsol \times
        \dualsol}) + \nabla \Vclf(x,z) w}{1+\lambda \circ
      \tilde{\alpha}_1^{-1}(\Vclf(x,z))}, \nonumber \\ &\le
    \scalebox{1.15}{$\frac{-\tilde{\alpha}_3(\NormTwo{(x,z)}_{\primalsol
          \times \dualsol})}{1 + \lambda \circ \tilde{\alpha}_1^{-1}
        \circ \tilde{\alpha}_2(\NormTwo{(x,z)}_{\pdsol})} +
      \frac{\lambda(\NormTwo{(x,z)}_{\primalsol \times \dualsol})}{1 +
        \lambda (\NormTwo{(x,z)}_{\primalsol \times \dualsol})}
      \NormTwo{w}$} \le -\rho(\NormTwo{(x,z)}_{\pdsol}) +
    \NormTwo{w}, \label{eq:W1_2}
  \end{align}
  where $\rho$ is the positive definite function given by
  \begin{align*}
    \rho(r) = \tilde{\alpha}_3(r)/(1+\lambda \circ
      \tilde{\alpha}_1^{-1} \circ \tilde{\alpha}_2(r)).
  \end{align*}  
  and we have used the fact that $\tilde{\alpha}_1^{-1}$ and
  $\tilde{\alpha}_2$ are positive definite. We now show that $V =
  \Vcircsep + \Vclfsep$ is an iISS Lyapunov function
  for~\eqref{eq:flow_robust} with respect
  to~$\pdsol$. First,~\eqref{eq:proper} is satisfied because $V$ is
  positive definite and radially unbounded with respect to $\pdsol$
  since (i) $\Vcircsep$ is positive definite and radially unbounded
  with respect to $\primalsol \times \dualsol$ and (ii) $\Vclfsep$ is
  positive definite with respect to $\pdsol$.
  Second,~\eqref{eq:dissipate} is satisfied as a result of the
  combination of~\eqref{eq:W1_1} and~\eqref{eq:W1_2}.
  Since $V$ satisfies the conditions of
  Theorem~\ref{th:iISS_lyap},~\eqref{eq:flow_robust} is iISS.
\end{IEEEproof}

Based on the discussion in Section~\ref{sec:not-ISS}, the iISS
property of~\eqref{eq:flow_robust} is an accurate representation of
the robustness of the dynamics, not a limitation of our analysis. A
consequence of iISS is that the asymptotic convergence of the dynamics
is preserved under finite energy disturbances~\cite[Proposition
  6]{EDS:98a}. In the case of~\eqref{eq:flow_robust}, a stronger
convergence property is true under finite variation disturbances
(which do not have finite energy).
The following formalizes this fact.

\begin{corollary}\longthmtitle{Finite variation
    disturbances}\label{cor:variation} 
  Suppose $w : \realsnonnegative \rightarrow \reals^{n} \times
  \reals^{m}$ is such that $\int_0^{\infty}\NormTwo{w(s) -
    \overline{w}} ds < \infty$ for some $\overline{w} =
  (\overline{w}_x,\overline{w}_z) \in \reals^{n}\times\reals^{m}$.
  Assume that $\primalsol(\overline{w}) \times \dualsol(\overline{w})$
  is nonempty and compact. Then each trajectory
  of~\eqref{eq:flow_robust} under the disturbance $w$ converges
  asymptotically to a point in $\primalsol(\overline{w}) \times
  \dualsol(\overline{w})$.

\end{corollary}
\begin{IEEEproof}
  Let $f_{\operatorname{dis,pert}}^v$ be the discontinuous
  saddle-point dynamics derived for the perturbed
  program~\eqref{eq:pert_prog} associated to $\overline{w}$ with
  additive disturbance $v: \realsnonnegative \rightarrow \reals^n
  \times \reals^m$.  By Corollary~\ref{cor:pert_prog},
  $\primalsol(\overline{w}) \times \dualsol(\overline{w}) \not=
  \emptyset$ is globally asymptotically stable for
  $f_{\operatorname{dis,pert}}^0$. Additionally, by
  Theorem~\ref{th:saddle_iISS} and since $\primalsol(\overline{w})
  \times \dualsol(\overline{w})$ is compact,
  $f_{\operatorname{dis,pert}}^v$ is iISS.  As a consequence,
  by~\cite[Proposition 6]{EDS:98a}, each trajectory of
  $f_{\operatorname{dis,pert}}^v$ converges asymptotically to a point
  in $\primalsol(\overline{w}) \times \dualsol(\overline{w})$ if
  $\int_0^{\infty}\NormTwo{v(s)} ds < \infty$.  The result now follows
  by noting that $\sdlw$ with disturbance $w$ is exactly
  $f_{\operatorname{dis,pert}}^v$ with disturbance $v = w -
  \overline{w}$ and that, by assumption, the latter disturbance
  satisfies $\int_0^{\infty}\NormTwo{v(s)} ds < \infty$.
\end{IEEEproof}

\myclearpage
\section{Robustness in recurrently connected
  graphs}\label{sec:connected}

In this section, we build on the iISS properties of the saddle-point
dynamics~\eqref{eq:flow} to study its convergence under communication
link failures. As such, agents do not receive updated state
information from their neighbors at all times and use the last known
value of their state to implement the dynamics.  The link failure
model we considered is described by recurrently connected graphs
(RCG), in which periods of communication loss are followed by periods
of connectivity. We formalize this notion next.

\begin{definition}\longthmtitle{Recurrently connected
    graphs}\label{def:connect} 
  Given a strictly increasing sequence of times
  $\{t_k\}_{k=0}^{\infty} \subset \reals_{\ge 0}$ and a base graph
  $\graph_b = (\vertices,\edges_b)$, we call $\graph(t) =
  (\vertices,\edges(t))$ \emph{recurrently connected with respect to
    $\graph_b$ and $\{t_k\}_{k=0}^{\infty}$} if $\edges(t) \subseteq
  \edges_b$ for all $t \in [t_{2k},t_{2k+1})$ while $\edges(t)
  \supseteq \edges_b$ for all $t \in [t_{2k+1},t_{2k+2})$, $k \in
  \nonnegint$.
\end{definition}

Intuitively, one may think of $\graph_b$ as a graph over
which~\eqref{eq:flow2} is distributed: during time intervals of the
form $[t_{2k},t_{2k+1})$, links are failing and hence the network
cannot execute the algorithm properly, whereas during time intervals
of the form $[t_{2k+1},t_{2k+2})$, enough communication links are
available to implement it correctly.
In what follows, and for simplicity of presentation, we only consider
the worst-case link failure scenario: i.e., if a link fails during the
time interval $[t_{2k},t_{2k+1})$, it remains down during its entire
duration.  The results stated here also apply to the general scenarios
where edges may fail and reconnect multiple times within a time
interval.

In the presence of link failures, the implementation of the evolution
of the $z$ variables, cf.~\eqref{eq:flow_z_ell}, across different
agents would yield in general different outcomes (given that different
agents have access to different information at different times). To
avoid this problem, we assume that, for each $\ell \in \{1,\dots,m\}$,
the agent with minimum identifier index,
\begin{align*}
  j = \sharedagent(\ell) := \min \{ i \in \{1,\dots,n\} : a_{\ell,i}
  \not= 0 \},
\end{align*}
implements the $z_{\ell}$-dynamics and communicates this value when
communication is available to its neighbors. Incidentally, only
neighbors of $j = \sharedagent(\ell)$ need to know $z_{\ell}$. With
this convention in place, we may describe the network dynamics under
link failures. Let $\failededges(k)$ be the set of failing
communication edges for $t \in [t_k,t_{k+1})$. In other words, if
$(i,j) \in \failededges(k)$ then agents $i$ and $j$ do not receive
updated state information from each other during the whole interval
$[t_k,t_{k+1})$. The nominal flow function of $i$ on a RCG for $t \in  
[t_k,t_{k+1})$ is
\begin{subequations}\label{eq:flow2_async}
  \begin{align}
    \anomflow_i&(x,z) = - c_i - \hspace{-4.7mm} \sum_{\substack{\ell =
        1
        \\
        (i,\sharedagent(\ell)) \notin \failededges(k)}}^m
    \hspace{-4.7mm} a_{\ell,i} z_{\ell} - \hspace{-4.7mm}
    \sum_{\substack{\ell = 1 \\ (i,\sharedagent(\ell)) \in
        \failededges(k)}}^m \hspace{-4.7mm} a_{\ell,i} z_{\ell}(t_k)
    - \sum_{\ell = 1}^m a_{\ell,i} \Big[ \hspace{-3mm}
    \sum_{\substack{j = 1 \\ (i,j) \notin \failededges(k)}}^n
    \hspace{-3mm} a_{\ell,j}x_j + \hspace{-3mm} \sum_{\substack{ j = 1
        \\
        (i,j) \in \failededges(k)}}^n \hspace{-3mm} a_{\ell,j}x_j(t_k)
    - b_\ell\Big]. \nonumber
  \end{align}
  Thus the $x_i$-dynamics during $[t_k,t_{k+1})$ for $i
  \in~\{1,\dots,n\}$~is
  \begin{align}
    \dot{x}_i &=
    \begin{cases}
      \anomflow_i(x,z), & \text{if } x_i > 0,
      \\
      \max\{0,\anomflow_i(x,z)\}, & \text{if } x_i = 0.
    \end{cases}
    \label{eq:flow2_x_async}
  \end{align}
  Likewise, the $z$-dynamics for $\ell \in \{1,\dots,m\}$ is
  \begin{align}
    \dot{z}_{\ell} = \hspace{-4.7mm} \sum_{\substack{i = 1
        \\
        (i,\sharedagent(\ell)) \notin \failededges(k)}}^n
    \hspace{-4.7mm} a_{\ell,i}x_i + \hspace{-4.7mm} \sum_{\substack{i
        = 1
        \\
        (i,\sharedagent(\ell)) \in \failededges(k)}}^n \hspace{-4.7mm}
    a_{\ell,i}x_i(t_k) - b_{\ell}.
  \end{align}
\end{subequations}
It is worth noting that~\eqref{eq:flow2_async} and~\eqref{eq:flow2}
coincide when $\failededges(k) = \emptyset$.  The next result shows
that the discontinuous saddle-point dynamics still converge under
recurrently connected graphs.

\begin{proposition}\longthmtitle{Convergence of saddle-point dynamics
    under RCGs}\label{prop:RCG}
  Let $\graph(t) = (\vertices,\edges(t))$ be recurrently connected
  with respect to $\graph_b = (\vertices,\edges_b)$ and
  $\{t_k\}_{k=0}^{\infty}$. Suppose that~\eqref{eq:flow2_async} is
  distributed over $\graph_b$ and $\Tdiscon := \sup_{k \in \nonnegint}
  (t_{2k+1}-t_{2k}) < \infty$. Let $t \mapsto (x(t),z(t))$ be a
  trajectory of~\eqref{eq:flow2_async}. Then there exists $\Tcon > 0$
  (depending on $\Tdiscon$, $x(t_0)$, and $z(t_0)$) such that $\inf_{k
    \in \nonnegint} (t_{2k+2}-t_{2k+1}) > \Tcon$ implies that
  $\NormTwo{(x(t_{2k}),z(t_{2k}))}_{\pdsol} \rightarrow 0$ as $k
  \rightarrow \infty$.
\end{proposition}
\begin{IEEEproof}
  The proof method is to (i) show that trajectories
  of~\eqref{eq:flow2_async} do not escape in finite time and (ii) use
  a $\mathcal{KL}$ characterization of asymptotically stable
  dynamics~\cite{CC-ART-RG:08} to find $\Tcon$ for which
  $\NormTwo{(x(t_{2k}),z(t_{2k}))}_{\pdsol} \rightarrow 0$ as $k
  \rightarrow \infty$. To prove (i), note that~\eqref{eq:flow2_async}
  represents a switched system of affine differential equations. The
  modes are defined by all $\kappa$-combinations of link failures (for
  $\kappa = 1,\dots,|\edges_b|$) and all $\kappa$-combinations of
  agents (for $\kappa = 1,\dots,n$). Thus, the number of modes is $d
  := 2^{|\edges_b|+n}$. Assign to each mode a number in the set
  $\{1,\dots,d\}$. Then, for any given $t \in [t_k,t_{k+1})$, the
  dynamics~\eqref{eq:flow2_async} is equivalently represented as
  \begin{align*}
    \scalebox{0.8}{$\left[
        \begin{array}{c}
          \dot{x}
          \\ 
          \dot{z}
        \end{array}
      \right]$}= P_{\sigma(t)}
    \scalebox{0.8}{$\left[
        \begin{array}{c}
          x 
          \\
          z 
        \end{array}
      \right]$} + q_{\sigma(t)}(x(t_k),z(t_k)),
  \end{align*}
  where $\sigma: \reals_{\ge 0} \rightarrow \{1,\dots,d\}$ is a
  switching law and $P_{\sigma(t)}$ (resp. $q_{\sigma(t)}$) is the
  flow matrix (resp. drift vector) of~\eqref{eq:flow2_async} for mode
  $\sigma(t)$. Let $\rho = \NormTwo{(x(t_0),z(t_0))}_{\primalsol
    \times \dualsol}$ and define
  \begin{align*}
    \tilde{q} := \hspace{-3mm} \max_{\substack{p \in \{1,\dots,d\}
        \\ \NormTwo{(x,z)}_{\pdsol} \le \rho}} \NormTwo{q_p(x,z)},
    \quad \text{and} \quad \tilde{\mu} := \hspace{-2mm} \max_{p \in
      \{1,\dots,d\}}\mu(P_p),
  \end{align*}
  where $\mu(P_p) = \lim_{h \rightarrow 0^+}
  \frac{\NormTwo{I-hP_p}^{-1}}{h}$ is the logarithmic norm of $P_p$.
  Both $\tilde{q}$ and $\tilde{\mu}$ are finite. Consider an arbitrary
  interval $[t_{2k},t_{2k+1})$ where
    $\NormTwo{(x(t_{2k}),z(t_{2k}))}_{\primalsol \times \dualsol} \le
    \rho$. In what follows, we make use of the fact that the
    trajectory of an affine differential equation $\dot{y} =
    \mathcal{A}y+\beta$ for $t \ge t_0$~is
  \begin{align}\label{eq:affine}
    y(t) &= e^{\mathcal{A}(t-t_0)}y(t_0) +
    \scalebox{1.2}{$\int_{t_0}^{t}$}e^{\mathcal{A}(t-s)}\beta ds.
  \end{align}
  Applying~\eqref{eq:affine}, we derive the following bound,
  \begin{align*}
    \NormTwo{(x(t_{2k+1}),z(t_{2k+1})) &- (x(t_{2k}),z(t_{2k}))}
    \\ &\le
    \NormTwo{(x(t_{2k}),z(t_{2k}))}(e^{\tilde{\mu}(t_{2k+1}-t_{2k})}-1)
     +
    \scalebox{1.2}{$\int_{t_{2k}}^{t_{2k+1}}$}e^{\tilde{\mu}(t_{2k+1}-s)}\tilde{q}ds,
    \\
    & \le (\rho + \tilde{q}/\tilde{\mu})(e^{\tilde{\mu}\Tdiscon}-1) =: M.
  \end{align*}
  In words, $M$ bounds the distance that trajectories travel on
  intervals of link failures. Also, $M$ is valid for all such
  intervals where $\NormTwo{(x(t_{2k}),z(t_{2k}))}_{\pdsol} \le
  \rho$. Next, we address the proof of (ii) by designing $\Tcon$ to
  enforce this condition. By definition,
  $\NormTwo{(x(t_{0}),z(t_{0}))}_{\pdsol} = \rho$. Thus,
  $\NormTwo{(x(t_{1}),z(t_{1})) - (x(t_{0}),z(t_{0}))} \le M$. Given
  that $\pdsol$ is globally asymptotically stable
  for~\eqref{eq:flow2_async} if $\failededges(k) = \emptyset$ (cf.
  Theorem~\ref{th:saddle_iISS}), \cite[Theorem
  3.13]{CC-ART-RG:08} implies the existence of $\beta \in
  \mathcal{KL}$ such that
  \begin{align*}
    \NormTwo{(x(t),z(t))}_{\pdsol} \le
    \beta(\NormTwo{(x(t_0),z(t_0))}_{\pdsol},t) .
  \end{align*}
  By~\cite[Proposition 7]{EDS:98a}, there exist $\theta_1,\theta_2 \in
  \mathcal{K}_{\infty}$ such that $ \beta(s,t) \le
  \theta_1(\theta_2(s)e^{-t})$. Thus,
  \begin{align*}
    \alpha(\NormTwo{(x(t_{2}),z(t_{2}))}_{\primalsol \times
      \dualsol})
     \le \theta_1(\theta_2(\NormTwo{(x(t_{1}),z(t_{1}))}_{\primalsol
      \times \dualsol})e^{-t_2+t_1}) \le \theta_1(\theta_2(\rho+M)e^{-t_2+t_1}).
  \end{align*}
  Consequently, if
  \begin{align*}
    t_2-t_1 > \Tcon := \ln \bigg(
    \frac{\theta_2(\rho+M)}{\theta_1^{-1}(\alpha(\rho))} \bigg) > 0,
  \end{align*}
  then $\NormTwo{(x(t_{2}),z(t_{2}))}_{\pdsol} < \rho$. Repeating this
  analysis reveals that $\NormTwo{(x(t_{2k+2}),z(t_{2k+2}))}_{\pdsol}
  < \NormTwo{(x(t_{2k}),z(t_{2k}))}_{\pdsol}$ for all $k \in
  \nonnegint$ when $t_{2k+2} - t_{2k+1} > \Tcon$. Thus
  \\ $\NormTwo{(x(t_{2k}),z(t_{2k}))}_{\pdsol} \rightarrow 0$ as $k
  \rightarrow \infty$ as claimed.
\end{IEEEproof}

\begin{remark}\longthmtitle{More general link
    failures}\label{rem:lack-robustness-jointly-connected} {\rm
    Proposition~\ref{prop:RCG} shows that, as long as the
    communication graph is connected with respect to $A$ for a
    sufficiently long time after periods of failure, the discontinuous
    saddle-point dynamics converge. We have observed in simulations,
    however, that the dynamics is not robust to more general link
    failures such as when the communication graph is never connected
    with respect to $A$ but its union over time is. We believe the
    reason is the lack of consistency in the $z-$dynamics for all time
    across agents in this case.  \oprocend }
\end{remark}

\section{Simulations}\label{sec:simulations}

Here we illustrate the convergence and robustness properties of the
discontinuous saddle-point dynamics.
We consider a finite-horizon optimal control problem for a network of
agents with coupled dynamics and underactuation. The network-wide
dynamics is open-loop unstable and the aim of the agents is to find a
control to minimize the actuation effort and ensure the network state
remains small.  To achieve this goal, the agents use the discontinuous
saddle-point dynamics~\eqref{eq:flow2}.  Formally, consider the
 finite-horizon optimal control problem,
\begin{subequations}\label{eq:ftoc}
  \begin{alignat}{2}
    &\min && \quad \sum_{\tau = 0}^T \NormTwo{x(\tau+1)}_1 +
    \NormTwo{u(\tau)}_1 \\ &\hspace{1.5mm}\text{s.t.} && \quad
    x(\tau+1) = Gx(\tau) + Hu(\tau), \quad \tau = 0, \dots
    T, \label{eq:ftoc2-network-evol} 
  \end{alignat}
\end{subequations}
where $x(\tau) \in \reals^N$ and $u(\tau) \in \reals^N$ is the network
state and control, respectively, at time $\tau$. The initial point
$x_i(0)$ is known to agent $i$ and its neighbors. The matrices $G \in
\reals^{N \times N}$ and $H = \diag(h) \in \reals^{N \times N}$, $h
\in \reals^N$, define the network evolution, and the network topology
is encoded in the sparsity structure of $G$.  We interpret each agent
as a subsystem whose dynamics is influenced by the states of
neighboring agents.  An agent knows the dynamics of its own subsystem
and its neighbor's subsystem, but does not know the entire network
dynamics.  A solution to~\eqref{eq:ftoc} is a time history of optimal
controls $(u_*(0),\dots,u_*(T)) \in (\reals^{N})^T$.

\begin{figure}[hbt!]
  \centering
  \subfigure[Network dynamics]{\includegraphics[trim=0 -15 0 0, width=.65\linewidth]{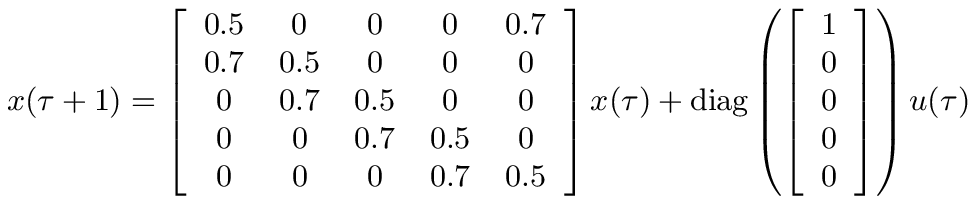}} \qquad 
  \subfigure[Communication topology]{\includegraphics[trim=20 40 10
    10, width=.23\linewidth]{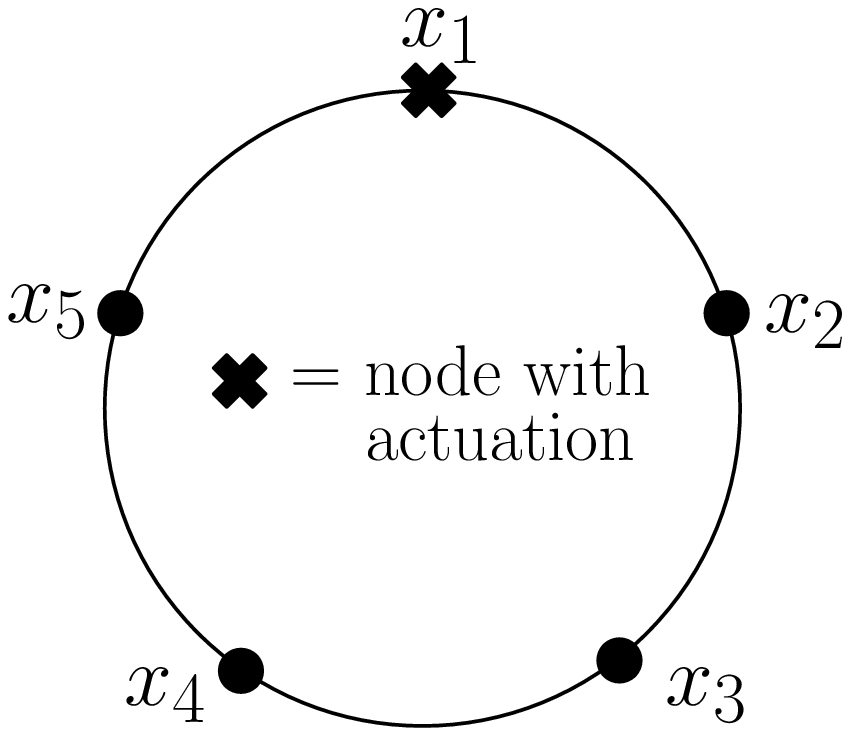}} 
  \caption{Network dynamics and communication topology of the
    multi-agent system. The network dynamics is underactuated and
    open-loop unstable but controllable. The presence of a
    communication link in (b) among every pair of agents whose
    dynamics are coupled in (a) ensures that the
    algorithm~\eqref{eq:flow2} is distributed over the communication
    graph.  }\label{fig:sim_dyn}
  \vspace*{-2ex}
\end{figure}

To express this problem in standard linear programming
form~\eqref{eq:standard_form}, we split the states into their positive
and negative components, $x(\tau) = x^+(\tau) - x^-(\tau)$, with
$x^+(\tau),x^-(\tau) \ge 0$ (and similarly for the inputs
$u(\tau)$). Then,~\eqref{eq:ftoc} can be equivalently formulated as
the following linear program,
\vspace{-4mm}
\begin{subequations}\label{eq:ftoc2}
  \begin{alignat}{2}
    &\min && \quad \sum_{\tau = 0}^T \sum_{i=1}^N x^+_i(\tau+1) +
    x^-_i(\tau+1) + u^+_i(\tau) + u^-_i(\tau)
    \\
    &\hspace{1.5mm}\text{s.t.}  && \quad x^+(\tau+1) - x^-(\tau) =
    G(x^+(\tau) - x^-(\tau)) + H(u^+(\tau) - u^-(\tau)), \quad \tau =
    0, \dots, T \\ & && \quad x^+(\tau+1), x^-(\tau+1), u^+(\tau),
    u^-(\tau) \ge 0, \quad \tau = 0, \dots, T
  \end{alignat}
\end{subequations}
The optimal control for~\eqref{eq:ftoc} at time $\tau$ is then
$u_*(\tau) = u^+_*(\tau) - u^-_*(\tau)$, where the vector
$(u^+_*(0),u^-_*(0),\dots,u^+_*(T),u^-_*(T))$ is a solution
to~\eqref{eq:ftoc2}, cf.~\cite[Lemma 6.1]{GBD:97}.

We implement the discontinuous saddle-point dynamics~\eqref{eq:flow2}
for problem~\eqref{eq:ftoc2} over the network of $5$ agents described
in Figure~\ref{fig:sim_dyn}.  To implement the
dynamics~\eqref{eq:flow2}, neighboring agents must exchange their
state information with each other.  In this example, each agent is
responsible for $2(T+1) = 24$ variables, which is independent of the
network size. This is in contrast to consensus-based distributed
optimization algorithms, where each agent would be responsible for
$2N(T+1) = 120$ variables, which grows linearly with the network
size~$N$. For simulation purposes, we implement the dynamics as a
single program in MATLAB\textsuperscript{\textregistered}, using a
first-order (Euler) approximation of the differential equation with a
stepsize of $0.01$. The CPU time for the simulation is $3.1824s$ on a
64-bit 3GHz Intel\textsuperscript{\textregistered}
Core\textsuperscript{TM} i7-3540M processor with 16GB of installed
RAM.
\begin{figure}[hbt!]
  \centering
  \subfigure[Computing the optimal control (with
  noise)]{\includegraphics[width=.495\linewidth]{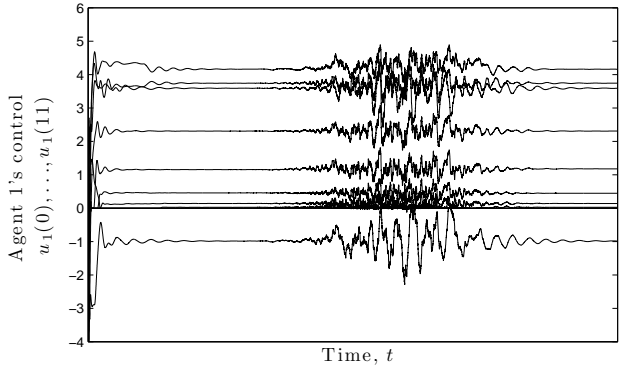}}
  \subfigure[Finite energy noise used in (a)]{\includegraphics[width=.495\linewidth]{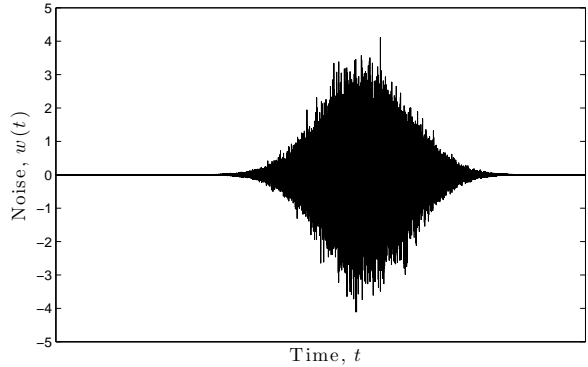}}
  \\
  \subfigure[Equality constraint violation in
  (a)]{\includegraphics[width=.495\linewidth]{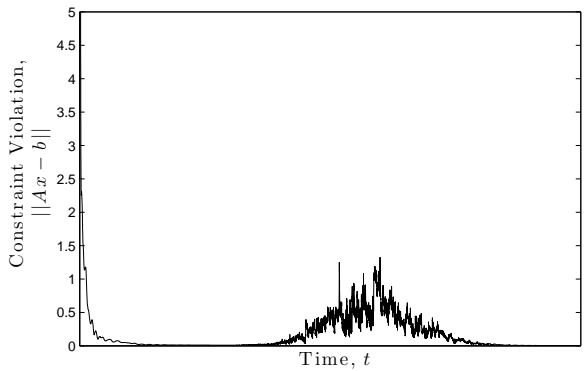}}
  \subfigure[Network evolution under optimal control found in
  (a)]{\includegraphics[width=.495\linewidth]{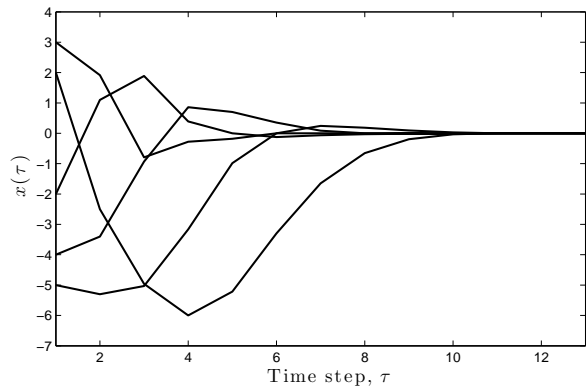}}
  \caption{Plot (a) shows the trajectories of the discontinuous
    saddle-point dynamics~\eqref{eq:flow_robust} subject to the noise
    depicted in (b) for agent $1$ as it computes its time history of
    optimal controls.  Plot (c) shows the associated equality
    constraint violation.  The asymptotic convergence of the
    trajectories appears to be exponential. The time horizon of the
    optimal control problem~\eqref{eq:ftoc2} is $T=11$. The $12$
    trajectories in (a) and (b) represent agent $1$'s evolving
    estimates of the optimal controls $u_1(0) ,\dots, u_1(11)$. The
    steady-state values achieved by these trajectories correspond to
    the solution of~\eqref{eq:ftoc}.  Once determined, these controls
    are then implemented by agent 1 and result in the network
    evolution depicted in (d).  The dynamics is initialized to a
    random point.  }\label{fig:sims1}
  \vspace*{-2ex}
\end{figure}
\begin{figure}[hbt!]
  \centering
  \includegraphics[width=.5\linewidth]{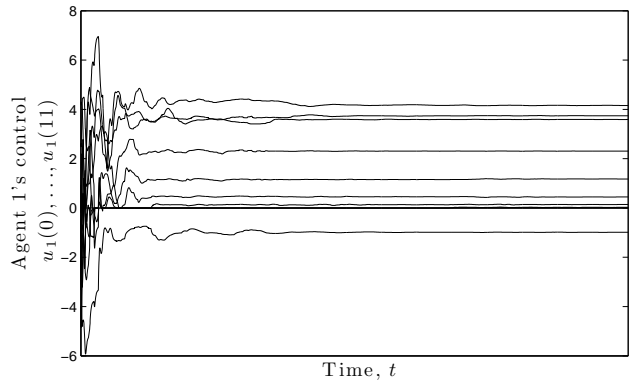}
    \caption{The trajectories of the discontinuous saddle-point
      dynamics~\eqref{eq:flow_robust} under a recurrently connected
      communication graph where a random number of random links failed
      during periods of disconnection. The simulation parameters are
      the same as in Figure~\ref{fig:sims1}. }\label{fig:sims2}
  \vspace*{-2ex}
\end{figure}

Note that, when implementing this dynamics, agent $i \in \until{5}$
computes the time history of its optimal control,
$u^-_i(0),u^+_i(0),\dots,u^-_i(T),u^+_i(T)$, as well as the time
history of its states,
$x^-_i(1),x^+_i(1),\dots,x^-_i(T+1),x^+_i(T+1)$. With respect to the
solution of the optimal control problem, the time history of states
are auxiliary variables used in the discontinuous dynamics and can be
discarded after the control is determined.  Figure~\ref{fig:sims1}
shows the results of the implementation of~\eqref{eq:flow2} when a
finite energy noise signal disturbs the agents'
execution. Clearly~\eqref{eq:flow2} achieves convergence initially in
the absence of noise. Then, the finite energy noise signal in
Figure~\ref{fig:sims1}(b) enters each agents' dynamics and disrupts
this convergence, albeit not significantly due to the iISS property
of~\eqref{eq:flow_robust} characterized in
Theorem~\ref{th:saddle_iISS}. Once the noise disappears, convergence
ensues. The constraint violation is plotted in
Figure~\ref{fig:sims1}(c). Once the time history of optimal controls
has been computed (corresponding to the steady-state values in
Figure~\ref{fig:sims1}(a)), agent 1 implements it, and the resulting
network evolution is displayed in Figure~\ref{fig:sims1}(d).  Agent 1
is able to drive the system state to zero, despite it being open-loop
unstable. Figure~\ref{fig:sims2} shows the results of implementation
in a recurrently connected communication graph and~\eqref{eq:flow2}
still achieves convergence as characterized in
Proposition~\ref{prop:RCG}.  The link failure model here is a random
number of random links failing during times of disconnection. The
graph is repeatedly connected for $1s$ and then disconnected for $4s$
(i.e., the ratio $\Tdiscon : \Tcon$ is $4:1$). The fact that
convergence is still achieved under this unfavorable ratio highlights
the strong robustness properties of the algorithm.

\myclearpage

\section{Conclusions}\label{sec:conclusions}

We have considered a network of agents whose objective is to have the
aggregate of their states converge to a solution of a general linear
program. We proposed an equivalent formulation of this problem in
terms of finding the saddle points of a modified Lagrangian
function. To make an exact correspondence between the solutions of the
linear program and saddle points of the Lagrangian we incorporate a
nonsmooth penalty term. This formulation has naturally led us to study
the associated saddle-point dynamics, for which we established the
point-wise convergence to the set of solutions of the linear program.
Based on this analysis, we introduced an alternative algorithmic
solution with the same asymptotic convergence properties. This
dynamics is amenable to distributed implementation over a multi-agent
system, where each individual controls its own component of the
solution vector and shares its value with its neighbors. We also
studied the robustness against disturbances and link failures of this
dynamics. We showed that it is integral-input-to-state stable but not
input-to-state stable (and, in fact, no algorithmic solution for
linear programming is). These results have allowed us to formally
establish the resilience of our distributed dynamics to disturbances
of finite variation and recurrently disconnected communication graphs.
Future work will include the study of the convergence rate of the
dynamics and its robustness properties under more general link
failures, the synthesis of continuous-time computation models with
opportunistic discrete-time communication among agents, and the
extension of our design to other convex optimization problems.  We
also plan to explore the benefits of the proposed distributed dynamics
in a number of engineering scenarios, including the smart grid and
power distribution, bargaining and matching in networks, and model
predictive control.


\bibliographystyle{ieeetr}%
\bibliography{alias,Main,Main-add,JC}

\begin{thebibliography}{10}

\bibitem{RA-RC-AC-LS:12}
R.~Alberton, R.~Carli, A.~Cenedese, and L.~Schenato, ``Multi-agent perimeter
  patrolling subject to mobility constraints,'' in {\em {A}merican {C}ontrol
  {C}onference}, (Montreal), pp.~4498--4503, 2012.

\bibitem{DPB:98}
D.~P. Bertsekas, {\em {N}etwork {O}ptimization: {C}ontinuous and {D}iscrete
  {M}odels}.
\newblock {A}thena {S}cientific, 1998.

\bibitem{MJ-SA-ME:06}
M.~Ji, S.~Azuma, and M.~Egerstedt, ``Role-assignment in multi-agent
  coordination,'' {\em International Journal of Assistive Robotics and
  Mechatronics}, vol.~7, no.~1, pp.~32--40, 2006.

\bibitem{BWC-AB-FD-BK-KR-GS-RB-MB-FA:12}
B.~W. Carabelli, A.~Benzing, F.~D\"urr, B.~Koldehofe, K.~Rothermel, G.~Seyboth,
  R.~Blind, M.~Burger, and F.~Allgower, ``Exact convex formulations of
  network-oriented optimal operator placement,'' in {\em {IEEE} Conf.\ on
  Decision and Control}, (Maui), pp.~3777--3782, Dec. 2012.

\bibitem{DRK-JP:64}
D.~R. Kuehn and J.~Porter, ``The application of linear programming techniques
  in process control,'' {\em IEEE Transactions on Applications and Industry},
  vol.~83, no.~75, pp.~423--427, 1964.

\bibitem{JT-ZH-MJ:07}
J.~Trdlicka, Z.~Hanzalek, and M.~Johansson, ``Optimal flow routing in multi-hop
  sensor networks with real-time constraints through linear programming,'' in
  {\em IEEE Conf. on Emerging Tech. and Factory Auto.}, pp.~924--931, 2007.

\bibitem{WFS:67}
W.~F. Sharpe, ``A linear programming algorithm for mutual fund portfolio
  selection,'' {\em Management Science}, vol.~13, no.~7, pp.~499--510, 1967.

\bibitem{GBD:63}
G.~B. Dantzig, {\em Linear Programming and Extensions}.
\newblock Princeton, NJ: Princeton University Press, 1963.

\bibitem{DB-JNT:97}
D.~Bertsimas and J.~N. Tsitsiklis, {\em Introduction to Linear Optimization},
  vol.~6 of {\em Optimization and Neural Computation}.
\newblock Belmont, MA: Athena Scientific, 1997.

\bibitem{SB-LV:09}
S.~Boyd and L.~Vandenberghe, {\em Convex Optimization}.
\newblock Cambridge University Press, 2009.

\bibitem{MB-GN-FB-FA:12}
M.~Burger, G.~Notarstefano, F.~Bullo, and F.~Allgower, ``A distributed simplex
  algorithm for degenerate linear programs and multi-agent assignment,'' {\em
  Automatica}, vol.~48, no.~9, pp.~2298--2304, 2012.

\bibitem{GN-FB:11}
G.~Notarstefano and F.~Bullo, ``Distributed abstract optimization via
  constraints consensus: Theory and applications,'' {\em IEEE Transactions on
  Automatic Control}, vol.~56, no.~10, pp.~2247--2261, 2011.

\bibitem{GY-RS:09}
G.~Yarmish and R.~Slyke, ``A distributed, scalable simplex method,'' {\em
  Journal of Supercomputing}, vol.~49, no.~3, pp.~373--381, 2009.

\bibitem{AN-AO-PAP:10}
A.~Nedic, A.~Ozdaglar, and P.~A. Parrilo, ``Constrained consensus and
  optimization in multi-agent networks,'' {\em IEEE Transactions on Automatic
  Control}, vol.~55, no.~4, pp.~922--938, 2010.

\bibitem{MZ-SM:12}
M.~Zhu and S.~Mart{\'\i}nez, ``On distributed convex optimization under
  inequality and equality constraints,'' {\em IEEE Transactions on Automatic
  Control}, vol.~57, no.~1, pp.~151--164, 2012.

\bibitem{JW-NE:11}
J.~Wang and N.~Elia, ``A control perspective for centralized and distributed
  convex optimization,'' in {\em {IEEE} Conf.\ on Decision and Control},
  (Orlando, Florida), pp.~3800--3805, 2011.

\bibitem{ROS-JAF-RMM:07}
R.~Olfati-Saber, J.~A. Fax, and R.~M. Murray, ``Consensus and cooperation in
  networked multi-agent systems,'' {\em Proceedings of the IEEE}, vol.~95,
  no.~1, pp.~215--233, 2007.

\bibitem{WR-RWB:08}
W.~Ren and R.~W. Beard, {\em Distributed Consensus in Multi-Vehicle Cooperative
  Control}.
\newblock Communications and Control Engineering, Springer, 2008.

\bibitem{FB-JC-SM:08cor}
F.~Bullo, J.~Cort\'es, and S.~Mart{\'\i}nez, {\em Distributed Control of
  Robotic Networks}.
\newblock Applied Mathematics Series, Princeton University Press, 2009.
\newblock Electronically available at http:/$\!$/coordinationbook.info.

\bibitem{MM-ME:10}
M.~Mesbahi and M.~Egerstedt, {\em Graph Theoretic Methods in Multiagent
  Networks}.
\newblock Applied Mathematics Series, Princeton University Press, 2010.

\bibitem{RC-GN:13}
R.~Carli and G.~Notarstefano, ``Distributed partition-based optimization via
  dual decomposition,'' in {\em {IEEE} Conf.\ on Decision and Control},
  (Firenze), pp.~2979--2984, Dec. 2013.

\bibitem{DPB-JNT:97}
D.~P. Bertsekas and J.~N. Tsitsiklis, {\em Parallel and Distributed
  Computation: Numerical Methods}.
\newblock Athena Scientific, 1997.

\bibitem{IN-JS:08}
I.~Necoara and J.~Suykens, ``Application of a smoothing technique to
  decomposition in convex optimization,'' {\em IEEE Transactions on Automatic
  Control}, vol.~53, no.~11, pp.~2674 -- 2679, 2008.

\bibitem{NG-JK:98}
N.~Garg and J.~Konemann, ``Faster and simpler algorithms for multicommodity
  flow and other fractional packing problems,'' in {\em Proceedings of the 39th
  Annual Symposium on Foundations of Computer Science}, (Palo Alto, CA),
  pp.~300--309, 1998.

\bibitem{DF-FP:10}
D.~Feijer and F.~Paganini, ``Stability of primal-dual gradient dynamics and
  applications to network optimization,'' {\em Automatica}, vol.~46,
  pp.~1974--1981, 2010.

\bibitem{KA-LH-HU:58}
K.~Arrow, L.~Hurwitz, and H.~Uzawa, {\em Studies in Linear and Non-Linear
  Programming}.
\newblock Stanford, California: Stanford University Press, 1958.

\bibitem{CC-ART-RG:08}
C.~Cai, A.~R. Teel, and R.~Goebel, ``Smooth {L}yapunov functions for hybrid
  systems part {II}: (pre)asymptotically stable compact sets,'' {\em IEEE
  Transactions on Automatic Control}, vol.~53, no.~3, pp.~734--748, 2008.

\bibitem{EDS:89-tac}
E.~D. Sontag, ``Further facts about input to state stabilization,'' {\em IEEE
  Transactions on Automatic Control}, vol.~35, pp.~473--476, 1989.

\bibitem{DA-EDS-YW:00}
D.~Angeli, E.~D. Sontag, and Y.~Wang, ``A characterization of integral
  input-to-state stability,'' {\em IEEE Transactions on Automatic Control},
  vol.~45, no.~6, pp.~1082--1097, 2000.

\bibitem{DB-DBB-CC:11}
D.~Bertsimas, D.~B. Brown, and C.~Caramanis, ``Theory and applications of
  robust optimization,'' {\em SIAM Review}, vol.~53, no.~3, pp.~464--501, 2011.

\bibitem{FHC:83}
F.~H. Clarke, {\em Optimization and Nonsmooth Analysis}.
\newblock Canadian Mathematical Society Series of Monographs and Advanced
  Texts, Wiley, 1983.

\bibitem{JC:08-csm-yo}
J.~Cort{\'e}s, ``Discontinuous dynamical systems - a tutorial on solutions,
  nonsmooth analysis, and stability,'' {\em {IEEE} Control Systems Magazine},
  vol.~28, no.~3, pp.~36--73, 2008.

\bibitem{DPB-AN-AEO:03}
D.~P. Bertsekas, A.~Nedi{\'c}, and A.~E. Ozdaglar, {\em Convex Analysis and
  Optimization}.
\newblock Belmont, MA: Athena Scientific, 1st~ed., 2003.

\bibitem{RD-PAS-RS:58}
R.~Dorfman, P.~A. Samuelson, and R.~Solow, {\em Linear programming in economic
  analysis}.
\newblock New York, Toronto, and London: McGraw Hill, 1958.

\bibitem{JPA-AC:84}
J.~P. Aubin and A.~Cellina, {\em Differential Inclusions}, vol.~264 of {\em
  Grundlehren der mathematischen Wissenschaften}.
\newblock New York: Springer, 1984.

\bibitem{WR:53}
W.~Rudin, {\em Principles of Mathematical Analysis}.
\newblock McGraw-Hill, 1953.

\bibitem{HJK-GGY:03}
H.~J. Kushner and G.~G. Yin, {\em Stochastic Approximation and Recursive
  Algorithms and Applications}, vol.~35 of {\em Applications of Mathematics:
  Stochastic Modelling and Applied Probability}.
\newblock New York: Springer, 2nd~ed., 2003.

\bibitem{RJBW:85}
R.~J.~B. Wets, ``On the continuity of the value of a linear program and of
  related polyhedral-valued multifunctions,'' {\em Mathematical Programming
  Study}, vol.~24, pp.~14--29, 1985.

\bibitem{EDS:98a}
E.~D. Sontag, ``Comments on integral variants of {ISS},'' {\em Systems \&
  Control Letters}, vol.~34, no.~1-2, pp.~93--100, 1998.

\bibitem{GBD:97}
G.~B. Dantzig, {\em Linear Programming: 1: Introduction}.
\newblock New York: Springer, 1997.

\bibitem{DPB:99}
D.~P. Bertsekas, {\em Nonlinear Programming}.
\newblock Belmont, MA: Athena Scientific, 2nd~ed., 1999.

\end{thebibliography}

\appendix

The following is a technical result used in the proof of
Theorem~\ref{th:no-ISS}.
\vspace{-2mm}
\begin{lemma}\longthmtitle{Property of feasible set}\label{lem:prop-feas}
  If $\{ Ax = b, x \ge 0 \}$ is non-empty and bounded then $\{A^Tz \ge
  0\}$ is unbounded.
\end{lemma}
\begin{IEEEproof}
  We start by proving that there exists an $\nu \in \reals^m$ such
  that $\{ Ax = b + \nu, x \ge 0 \}$ is empty. Define the vector $s
  \in \reals^n$ component-wise as $s_i = \max_{\{Ax = b, x \ge 0\}}
  x_i$.  Since $\{Ax = b, x \ge 0\}$ is compact and non-empty, $s$ is
  finite. Next, fix $\epsilon > 0$ and let $\nu = -A(s+\epsilon
  \mathbbm{1}_n)$. Note that $Ax = b + \nu$ corresponds to $A(x + s +
  \epsilon \mathbbm{1}_n) = b$, which is a shift by $s + \epsilon
  \mathbbm{1}_n$ in each component of $x$. By construction, $\{ Ax = b
  + \nu, x \ge 0 \}$ is empty. Then, the application of Farkas'
  Lemma~\cite[pp.~263]{SB-LV:09} yields that there exists $\hat{z} \in
  \reals^m$ such that $A^T\hat{z} \ge 0$ and $(b+\nu)^T\hat{z} < 0$
  (in particular, $(b+\nu)^T\hat{z} < 0$ implies that $\hat{z} \not=
  0$). For any $\lambda \in \reals_{\ge 0}$, it holds that
  $A^T(\lambda \hat{z}) \ge 0$, and thus $\lambda \hat{z} \in \{A^Tz
  \ge 0\}$, which implies the result.
\end{IEEEproof}

The proof of Theorem~\ref{th:saddle_iISS} makes use of the following
result from~\cite[Proposition B.25]{DPB:99}.
\vspace{-2mm}
\begin{theorem}\longthmtitle{Danskin's
    Theorem}\label{th:danskin} Let $Y \subset \reals^m$ be compact and
  convex. Given $g:\reals^n \times Y \rightarrow \reals$, suppose that
  $x \mapsto g(x,y)$ is differentiable for every $y \in Y$,
  $\partial_x g $ is continuous on $\reals^n \times Y$, and $y \mapsto
  g(x,y)$ is strictly convex and continuous for every $x \in
  \reals^n$. Define $f: \reals^n \rightarrow \reals$ by $f(x) =
  \min_{y \in Y} g(x,y)$.  Then, $\nabla f(x) = \partial_x g(x,y) |_{y
    = y_*(x)}$, where $y_*(x) = \argmin_{y \in Y} g(x,y)$.

\end{theorem}

\end{document}